\documentclass{amsbook}
\usepackage{amssymb}
\usepackage{amsfonts}

\makeatletter
\numberwithin{equation}{section}
\numberwithin{figure}{section}
\numberwithin{equation}{section}
\numberwithin{figure}{section}
\numberwithin{equation}{section}

\theoremstyle{definition}

\newtheorem*{acknowledgements*}{Acknowledgements}
\theoremstyle{remark}

\numberwithin{theorem}{section}   
\input{tcilatex}

\begin{document}
\frontmatter
\title[Ovchinnikov's Theorem]{Generalizing Ovchinnikov's Theorem}
\author{Jaak Peetre}
\address{Lund Institute of Technology, Lund, Sweden}
\author{Per G. Nilsson (Typist)}
\address{Nilsson: Stockholm, Sweden}
\email{pgn@plntx.com}
\date{\today }

\begin{abstract}
This note is an (exact) copy of the report of Jaak Peetre, "Generalizing
Ovchinnikov's Theorem". Published as Technical Report, Lund (1981). Some
more recent general references have been added, some references updated
though (in \textit{italics}) and some misprints corrected.
\end{abstract}

\maketitle

\begin{center}
\bigskip

{\Huge Generalizing Ovchinnikov's Theorem.}

{\Huge Jaak Peetre}
\end{center}

\section*{0. Introduction}

\pagenumbering{arabic}

Whereas in the past much work in interpolation spaces more or less
consciously has been connected with direct applications to various questions
in Analysis the last 4 - 5 years have witnessed a considerable progress on
the theoretical side. I am in the first place thinking of the results of
Ovchinnikov (including joint work with Dmitriev) and the important recent
paper by Janson - see the bibliography \cite{23-Ovc}, \cite{24-Ovc}, \cite%
{25-Ovc},\cite{9-DmKrOv} \cite{10-DmOv} and \cite{16-Jan} respectively.
Unfortunately most of Ovchinnikov's papers are hard to read and so is to
some extent Janson's too. A common basis for all this work is a fundamental
theorem by Aronszajn and Gagliardo \cite{1-AG65}, published already in 1965,
but which for quite a long time seemed to have fallen almost into oblivion%
\footnote{%
For references, post 1981, see notably Brudnyi-Krugljak \cite[Chapter 2,
Chapter 4]{33-BK}, Brudnyi-Krein-Semenov \cite{34-BrKrSe}, Kaijer-Pelletier 
\cite{35-KaPe},\cite{36-KaPe-II},, Ovchinnikov \cite{38-Ov84}, Nilsson \cite%
{37-Ni94} and the references listed there (especially \cite{34-BrKrSe}).}%
.Perhaps we are finally beginning to gain a better understanding of the true
inner meaning of the various interpolation methods which often have come up
in a rather ad hoc fashion\footnote{%
See Notes: \underline{Introductio\.{n}}:$\left\langle 1\right\rangle .$}. $%
\medskip $

\underline{0.1} A central result in Ovchinnikov's work is the following
interpolation theorem \cite{23-Ovc}. Let $T$ be any bounded linear operator
mapping the couple of weighted $l^{\infty }$ spaces $l^{\infty }\left( 
\overline{w}\right) =\left( l^{\infty }\left( w^{0}\right) ,l^{\infty
}\left( w^{1}\right) \right) $ into the couple of weighted $l^{1}$ spaces $%
l^{1}\left( \overline{w}\right) =\left( l^{1}\left( w^{0}\right)
,l^{1}\left( w^{1}\right) \right) $ in symbols: $T:l^{\infty }\left( 
\overline{w}\right) \rightarrow l^{1}\left( \overline{w}\right) .$ Then $T$
maps the space $l^{\infty }\left( w^{0}/\rho \left( w^{1}/w^{0}\right)
\right) $ into $l^{1}\left( w^{0}/\rho \left( w^{1}/\rho \left( w^{1}\right)
\right) \right) ,$ in symbols $T:l^{\infty }\left( w^{0}/\rho \left(
w^{1}/w^{0}\right) \right) \rightarrow l^{1}\left( w^{0}/\rho \left(
w^{1}/w^{0}\right) \right) .$ Here $\rho $ is any pseudo-concave function $%
\left( \rho \in \mathcal{P}\right) .$ Ovchinnikov's original proof was based
on Grothendieck's fundamental theorem \cite{14-Gr}. Subsequently Janson \cite%
{16-Jan} managed to give another more elementary proof which in particular
does not involve the latter result (cf. sec. \underline{5} of this paper).
(If $\rho \in \mathcal{P}^{+-}$ a proof based on Xincin's (= Khintchine's)
inequality was found by Gustavsson \cite{15-Gu}.) In connection with his
interpolation theorem Ovchinnikov \cite{23-Ovc} had introduced three new
interpolation methods, indexed by the letters $l,m$ and $u$ (l for "lower",
m for "middle" and u for "upper"), which fall in the general scheme of
Aronszajn and Gagliardo \cite{1-AG65}. Janson \cite{16-Jan} in the other
hand clarifies the role of l and u. However m seems now to fall out of the
picture. $\medskip $

\underline{0.2} The purpose of this paper is to extend Ovchinnikov's theorem
in various directions. In our first rather trivial generalization (sec. 
\underline{1}) we try to replace the spaces $l^{\infty }\left( w^{0}/\rho
\left( w^{1}/w^{0}\right) \right) $ and $l^{1}\left( w^{0}/\rho \left(
w^{1}/w^{0}\right) \right) $ by general interpolation (or more generally
only intermediate) spaces $A$ and $B$ with respect to $l^{\infty }\left( 
\overline{w}\right) $ and $l^{1}\left( \overline{w}\right) $ respectively.
Our result is roughly speaking that $T:l^{\infty }\left( \overline{w}\right)
\rightarrow l^{1}\left( \overline{w}\right) ,$ where $T$ is a bounded linear
operator implies $T:A\rightarrow B$ if and only if this assertion is already
true when $T$ is a multiplier transform. Next (sec. \underline{2}) we
consider vector valued and the continuous case. As a corollary (\underline{%
2.5} theorem) we get the following Ovchinnikov type theorem for Besov
couples: If $T:B_{p}^{\overline{s},\infty }\left( \mathbb{R}^{n}\right) 
\overset{\text{def}}{=}\left( B_{p}^{s_{0},\infty }\left( \mathbb{R}%
^{n}\right) ,B_{p}^{s_{1},\infty }\left( \mathbb{R}^{n}\right) \right)
\rightarrow B_{q}^{\overline{t},1}\left( \mathbb{R}^{m}\right) \overset{%
\text{def }}{=}\left( B_{q}^{t_{0},1}\left( \mathbb{R}^{m}\right)
,B_{q}^{t_{1},1}\left( \mathbb{R}^{m}\right) \right) $ then $%
T:B_{p}^{s,\infty }\left( \mathbb{R}^{n}\right) \rightarrow
B_{q}^{t,1}\left( \mathbb{R}^{m}\right) $ where $s=\left( 1-\theta \right)
s_{0}+\theta s_{1},t=\left( 1-\theta \right) t_{0}+\theta t_{1}$ with $%
\theta \in \left( 0,1\right) .$ In the following section (sec. \underline{3}%
) we introduce an abstract setting for Ovchinnikov type theorems. If $%
\overline{X}=\left( X_{0},X_{1}\right) $ and $\overline{Y}=\left(
Y_{0},Y_{1}\right) $ are general Banach couples we say they satisfy
condition $\left( O\right) $ if $T:\overline{X}\rightarrow \overline{Y}$
implies $T:K^{A}\left( \overline{X}\right) \rightarrow J^{B}\left( \overline{%
Y}\right) $ where $A$ and $B$ are any two of the spaces of sec. \underline{1}
and $K^{A}$ and $J^{B}$ stands for the corresponding $K$ and $J$ functors
respectively. We put this in relation with the work of Janson \cite{16-Jan}.
Janson has observed - in the special case $A=l^{\infty }\left( \rho \left( 
\overline{w}\right) \right) $ and $B=l^{1}\left( \rho \left( \overline{w}%
\right) \right) ,$which is really the only case which matters - that%
\footnote{%
See Notes: \underline{Introductio\.{n}}:$\left\langle 2\right\rangle .$} $%
K^{A}$ is the maximal functor with $K^{A}\left( l^{\infty }\left( \overline{w%
}\right) \right) \subseteq A$ and $J^{B}$ the minimal functor with $%
J^{B}\left( l^{1}\left( \overline{w}\right) \right) \supset B$ - maximal and
minimal of course in the sense of the Aronszajn-Gagliardo theorem \cite%
{1-AG65} . Let now $G^{A}$ be the minimal functor with $G^{A}\left(
l^{\infty }\left( \overline{w}\right) \right) \supset A$ and $H^{B}$ the
maximal functor with $H^{B}\left( l^{1}\left( \overline{w}\right) \right)
\subseteq B.$ We say that the couple $\overline{X}$ is of type $\left(
l\right) $ if $G^{A}\left( \overline{X}\right) =K^{A}\left( \overline{X}%
\right) $ for all $A$ and that the couple $\overline{Y}$ is of type $\left(
u\right) $ if $H^{B}\left( \overline{Y}\right) =J^{B}\left( \overline{Y}%
\right) $ for all $B.$ Then trivially, if $\overline{X}$ is of type $\left(
l\right) $ and $\overline{Y}$ of type $\left( u\right) ,$ the couples $%
\overline{X}$ and $\overline{Y}$ satisfies condition $\left( O\right) .$ For
instance, we show that Marcinkiewicz couples are of type $\left( l\right) $
and Lorentz couples of type $\left( u\right) .$Thus we get an Ovchinnikov
type theorem with Marcinkiewicz and Lorentz spaces instead of weighted $%
l^{\infty }$ and $l^{1}.$ In Sec. \underline{4} we relate those notions with
the idea of nuclearity. In this context we give also a direct proof of
Ovchinnikov's theorem \cite{23-Ovc}, along the lines of the original one,
thus relying to Grothendieck's fundamental theorem \cite{14-Gr}. En
revanche, in a attempt to make this paper as well-contained as possible, we
include what is essentially Janson's proof \cite{16-Jan} but free of the
context of Aronszajn-Gagliardo functors - whether this really is an
advantage from the conceptual point of view, we leave to the reader to
judge. Finally in the short sec. \underline{6} , pointing out some open
problems, we outline a program for further research.$\medskip $

\underline{0.3.} Some of our terminology is put forward in Sec. $\frac{1}{2}$
- obtained thus by interpolation between the present sec. \underline{0} and
sec. \underline{1}. As far as practicable we try to adhere to the one of
Bergh-L\"{o}fstr\"{o}m \cite{3-BeL0}.$\medskip $

We make the following convention: $\frac{a}{b}<\infty $ is interpreted to
mean $a=0$ if $b=0.\medskip $

The sign $\approx $ further stands as usual for equivalence. Thus $u\approx
v $ where $u$ and $v$ are positive quantities, means that we have
inequalities $u\leq C_{1}v,v\leq C_{2}u;$ what the constants $C_{1}$ and $%
C_{1}$ depend on is usually clear from the context.$\medskip $

\underline{0.4} The principal results of this paper - in a somewhat
premature form though - were already announced at the Scandinavian
Mathematical Congress held at Aarhus, Denmark in aug. 1980 \cite{26-Pe}. In
the meantime we have read Janson \cite{16-Jan}, which has lead to
considerable improvements.

\section*{$\frac{1}{2}$. Definitions.}

\underline{$\frac{1}{2}.1$}$.$ Let $\overline{X}=\left( X_{0},X_{1}\right) $
and $\overline{Y}=\left( Y_{0},Y_{1}\right) $ be two Banach pairs. As usual
we put $\Delta \overline{X}=X_{0}\cap X_{1}$ and $\Sigma \overline{X}%
=X_{0}+X_{1}$ (similarly for $\overline{Y}).$ Moreover $\Sigma _{0}\overline{%
X}$ denote the closure of $\Delta X$ in $\Sigma \overline{X}.$ If $\Sigma 
\overline{X}$ is \underline{regular} (i.e. if $\Delta \overline{X}$ is dense
in both $X_{0}$ and $X_{1}$) then $\Sigma _{0}\overline{X}=\Sigma \overline{X%
}$.$\medskip $

We say that we have a \underline{bounded linear} operator from $\overline{X}$
into $\overline{Y}$ if there is given a linear map $T$ from $\Sigma 
\overline{X}$ into $\Sigma \overline{Y}$ such that $T\left( X_{i}\right)
\subseteq Y_{i}$ ($i=0,1),$ the restriction $T_{i}=T\mid _{X_{i}}$of $T$ to $%
X_{i}$ defining a continuous linear operator from $X_{i}$ into $Y_{i}$ in
the ordinary sense. We then write $T:\overline{X}\rightarrow \overline{Y}$
or more completely $T:\overline{X}\overset{b}{\rightarrow }\overline{Y}.$ We
define the \underline{operator norm} of $T$ by $\left\Vert T\right\Vert
=\max \left( \left\Vert T_{0}\right\Vert ,\left\Vert T_{1}\right\Vert
\right) $ where $\left\Vert T_{i}\right\Vert $ is the usual operator norm of 
$T_{i},$ i.e. $\left\Vert T_{i}\right\Vert =\sup_{0\neq x\in
X_{i}}\left\Vert Tx\right\Vert _{Y_{i}}/\left\Vert x_{i}\right\Vert _{X_{i}}$
($i=0,1).\medskip $

We say that $T$ is \underline{nuclear} operator from $\overline{X}$ into $%
\overline{Y}$ and we use the notation $T:\overline{X}\overset{n}{\rightarrow 
}\overline{Y}$ , if $T$ is a linear map from $\Sigma \overline{X}$ into $%
\Sigma \overline{Y}$ such that there exists a family $\left\{ b_{n}\right\} $
in $\Delta \overline{Y}$ and a family $\left\{ l_{n}\right\} $ in $\left(
\Sigma \overline{X}\right) ^{^{\prime }}$ with\footnote{%
See Notes \underline{$\frac{1}{2}$}:$\left\langle 1\right\rangle .$} 
\begin{equation}
\left( 1\right) :\dsum \max \left\{ \left\Vert l_{n}\right\Vert
_{X_{0}^{^{\prime }}}\left\Vert b_{n}\right\Vert _{Y_{0}},\left\Vert
l_{n}\right\Vert _{X_{1}^{^{\prime }}}\left\Vert b_{n}\right\Vert
_{Y_{1}}\right\} <\infty   \tag{1}
\end{equation}%
such that $Ta=\dsum l_{n}\left( a\right) b_{n}$ for any $a\in \Sigma 
\overline{X}.$ (Here $n$ rungs through some index set $I,$which we without
loss of generality may take at most denumerable infinite.). Note that the
series $\dsum l_{n}\left( a\right) b_{n}$ is always summable in $\Sigma 
\overline{Y}.$ Indeed writing $\alpha =a_{0}+a_{1}$ with $a_{i}\in X_{i}$ $%
\left( i=0,1\right) $ we see that each of the series $\dsum l_{n}\left(
a_{i}\right) b_{n}$ is summable in $Y_{i},$ so the sum $Ta$ is really an
element of $\Sigma _{0}\overline{Y}.$ This also shows that $T:\overline{X}%
\overset{n}{\rightarrow }\overline{Y}$ implies $T:\overline{X}\overset{b}{%
\rightarrow }\overline{Y}.$ We define the \underline{\emph{nuclear norm}} of 
$T,$ denoted by $\left\Vert T\right\Vert _{n},$ as the infimum of all
expressions appearing in $\left( 1\right) .$ Clearly $\left\Vert
T\right\Vert \leq \left\Vert T\right\Vert _{n}.\medskip $

\underline{\emph{Remark.}}\emph{\ }If $T:\overline{X}\overset{n}{\rightarrow 
}\overline{Y}$ then already each of the maps $T:X_{i}\rightarrow Y_{i}$ $%
\left( i=0,1\right) \footnote{%
See Notes \underline{$\frac{1}{2}$}:$\left\langle 2\right\rangle .$}$ is
nuclear; we may say that $T$ is "separately nuclear". But the converse is
not true. There are counter-examples showing that a separately nuclear
operator need not to be nuclear (in our sense).$\medskip $

\underline{$\frac{1}{2}.2$}$.$ We introduce the following "orderings"
between elements of $\Sigma \overline{Y}$ and $\Sigma \overline{X}$ (or $%
\Sigma _{0}\overline{Y}$ and $\Sigma _{0}\overline{X}),$ $\lambda $ denoting
a fixed number $>1.$

$\left( i\right) :y<_{b}x$ means that for any $\epsilon >0$ there exists an
operator $T:\overline{X}\overset{b}{\rightarrow }\overline{Y}$ with $%
\left\Vert T\right\Vert <1+\epsilon $ such that $y=Tx.$

$\left( ii\right) :y<_{n}x$ means that for any $\epsilon >0$ there exists an
operator $T:\overline{X}\overset{n}{\rightarrow }\overline{Y}$ with $%
\left\Vert T\right\Vert _{n}<1+\epsilon $ such that $y=Tx$

$\left( iii\right) :y<_{K\left( \lambda \right) }x$ means that $K\left(
\lambda ^{k},y\right) \leq K\left( \lambda ^{k},x\right) $ for all $k\in 
\mathbb{Z}$.

$\left( iv\right) :y<_{J\left( \lambda \right) }x$ means that for every 
\underline{representation}\emph{\ }$\widehat{x}=\left\{ x_{k}\right\} _{k\in 
\mathbb{Z}}$ of $x\footnote{%
See Notes \underline{$\frac{1}{2}$}:$\left\langle 3\right\rangle .$}$ and
every $\epsilon >0$ there exists a representation $\widehat{y}=\left\{
y_{k}\right\} _{k\in \mathbb{Z}}$ of $y$ such that $J\left( \lambda
^{k},y_{k}\right) \leq \left( 1+\epsilon \right) J\left( \lambda
^{k},x_{k}\right) $ for all $k\in \mathbb{Z}$.

$\left( v\right) :y<_{J/K\left( \lambda \right) }x$ means that for every $%
\epsilon >0$ there exists a representation $\widehat{y}=\left\{
y_{k}\right\} _{k\in \mathbb{Z}}$ of $y$ such that $\dsum_{k\in \mathbb{Z}}%
\frac{J\left( \lambda ^{k},y_{k}\right) }{K\left( \lambda ^{k},x\right) }%
<1+\epsilon .$

Orderings seem to be destined to play a prominent role in interpolation
theory. E.g. \cite{6-CwPe} is mainly devoted to a study of the orderings $%
\left( iii\right) $ and $\left( iv\right) .$

The orderings $\left( i\right) ,\left( iii\right) $ and $\left( iv\right) $
are \underline{transitive} in the obvious sense. Let $\overline{Z}=\left(
Z_{0},Z_{1}\right) $ be another Banach pair (besides $\overline{X}$ and $%
\overline{Y}$). Then e.g. $z<_{b}y,y<_{b}x$ implies $z<_{b}x.$ The orderings 
$\left( ii\right) $ and $\left( v\right) $ have the \underline{ideal property%
} with respect to the ordering $\left( i\right) .$That is 
\begin{eqnarray*}
z &<&_{b}y,y<_{n}x\Longrightarrow z<_{n}x, \\
z &<&_{n}y,y<_{b}x\Longrightarrow z<_{n}x.
\end{eqnarray*}%
On the other hand the ordering $\left( v\right) $ has the \underline{left
ideal property} with respect to ordering $\left( iv\right) $ and the right
ideal property with respect to the ordering $\left( iii\right) .$That is;%
\begin{eqnarray*}
z &<&_{J\left( \lambda \right) }y,y<_{J/K\left( \lambda \right)
}x\Longrightarrow z<_{J/K\left( \lambda \right) }x; \\
z &<&_{J/K\left( \lambda \right) }y,y<_{K\left( \lambda \right)
}x\Longrightarrow z<_{J/K\left( \lambda \right) }x.
\end{eqnarray*}%
This is also true for the ordering $\left( ii\right) $ but this not quite
obvious. For the proof sec. \underline{4.1} cor. 2.$\medskip $

Let $<_{0}\left( o=b,n,...\right) $ be any of the orderings $\left( i\right)
-\left( v\right) .$Then we use the symbol $<<_{0}$ in the following sense. $%
y<<_{0}x$ means that $y<_{0}cx$ for some $c.$(Sometimes it is intended that
the constant $c$ depends only on the pairs involved. This will in most cases
be clear from the context.) In the case of the orderings, $\left( iii\right)
-\left( v\right) $ the value of $\lambda $ is then immaterial so we write
simply $y<<_{K}x$ in place of $y<<_{K\left( \lambda \right) }x$ etc.$%
\medskip $

\underline{$\frac{1}{2}.3$}$.$We recall now the definition of the
"classical" $K-$ and $J-$ spaces (see e.g. \cite{6-CwPe}; cf. \cite{3-BeL0}, 
\cite{4-BrKr}), in fact in a slightly more general form as suggested by the
work of Janson \cite{16-Jan}. Let $A$ be any intermediate space with respect
to the weighted $l^{\infty }$ couple $l^{\infty }\left( \overline{w}\right)
=\left( l^{\infty }\left( w^{0}\right) ,l^{\infty }\left( w^{1}\right)
\right) $ and $B$ be one with respect to the weighted $l^{1}$ couple $%
l^{1}\left( \overline{w}\right) =\left( l^{1}\left( w^{0}\right)
,l^{1}\left( w^{1}\right) \right) ,\overline{w}=\left( w^{0},w^{1}\right)
=\left( \left\{ w_{n}^{0}\right\} ,\left\{ w_{n}^{1}\right\} \right) $ be
any given couple of weight sequences.\footnote{%
See Notes \underline{$\frac{1}{2}$}:$\left\langle 4\right\rangle .$} Let $%
\overline{X}=\left( X_{0},X_{1}\right) $ be any Banach couple. Then we set%
\begin{equation*}
x\in K^{A}\left( \overline{X}\right) \text{ iff }\left\{ w_{n}^{0}K\left( 
\frac{w_{n}^{0}}{w_{n}^{1}},x\right) \right\} \in A,
\end{equation*}%
\begin{eqnarray*}
x &\in &J^{B}\left( \overline{Y}\right) \text{ iff for some representation }%
\widehat{x}=\left\{ x_{n}\right\} ~\text{of }x\text{ holds} \\
\left\{ w_{n}^{0}J\left( \frac{w_{n}^{0}}{w_{n}^{1}},x\right) \right\} &\in
&B.
\end{eqnarray*}%
The properties of these spaces are well-known. If we want to emphasize the
dependence on $\overline{w}$ we use the more complete notation $K_{\overline{%
w}}^{A}$ and $J_{\overline{w}}^{B},$in place of $K^{A}$ and $J^{B}$
respectively. (From Sec. 3 on we shall again specialize to the traditional
case $w_{n}^{0}=1,w_{n}^{1}=\lambda ^{-n}.)$. Sometimes we also put,
following the historical tradition, $K^{A}\left( \overline{X}\right) =%
\overline{X}_{A;K}$ and $J^{B}\left( \overline{Y}\right) =\overline{Y}%
_{B;J}. $ In particular let $A$ and/or $B$ be the space $l^{p}\left( \rho
\left( \overline{w}\right) \right) $ corresponding to the norm 
\begin{equation*}
\left\Vert c\right\Vert =\left( \dsum_{n}\left( w_{n}^{0}\left\vert
c_{n}\right\vert /\rho \left( \frac{w_{n}^{0}}{w_{n}^{1}}\right) \right)
^{p}\right) ^{1/p}\text{ (with }c=\left\{ c_{n}\right\} )
\end{equation*}%
where $p\in \left[ 1,\infty \right] $ (with the usual interpretation for $%
p=\infty $) and $\rho $ any positive function, usually pseudo-concave though
($\rho \in \mathcal{P)}$ . Then we write $\overline{X}_{\rho p;K}$ and $%
\overline{X}_{\rho p;J}$ respectively for these spaces. If $\rho \left(
t\right) =t^{\theta }$ $\left( \theta \in \left[ 0,1\right] \right) $ we get
the space $\overline{X}_{\theta p;K}$ and $\overline{X}_{\theta p;J}.$ (If $%
\theta \in \left( 0,1\right) $ we can, in view of the equivalence theorem,
omit the additional subscripts $K$ and $J$ writing $\overline{X}_{\theta p}$
for both of them.)$\medskip $

\underline{$\frac{1}{2}.4$}$.$ We mention also some basic facts connected
with the Aronszajn-Gargliardo theorem \cite{1-AG65} already referred to in
the introduction.$\medskip $

Let $\overline{A}=\left( A_{0},A_{1}\right) $ and $\overline{B}=\left(
B_{0},B_{1}\right) $ be any given two Banach couples and $A$ and $B$
intermediate spaces with respect to $\overline{A}$ and $\overline{B}$
respectively. Then there exists a minimal interpolation functor $Orb_{A}$
(or $Orb_{A,\overline{A}})$ such that $Orb_{A}\left( \overline{A}\right)
\supseteq A$ and a maximal one $Corb_{B}$ (or $Corb_{B,\overline{B}})$ such
that $Corb_{B}\left( \overline{B}\right) \subseteq B.$ If $\overline{X}%
=\left( X_{0},X_{1}\right) $ is a generic Banach couple we say that $%
Orb_{A}\left( \overline{X}\right) $ is the \underline{orbit} of $A$ in $%
\overline{X}$ and that $Corb_{B}\left( \overline{X}\right) $ is the 
\underline{coorbit} of $B$ in $\overline{X}.$ It is easy to see that if $A$
and $B$ are \underline{relative interpolation}\emph{\ }spaces with respect
to $\overline{A}$ and $\overline{B}$ then $Orb_{A}\subseteq
Corb_{B}.\medskip $

In particular if $\overline{A}=l^{\infty }\left( \overline{w}\right) $ or $%
B=l^{1}\left( \overline{w}\right) $ we put $G^{A}=Orb_{A}$ and $%
H^{B}=Corb_{B}$ respectively. In the special case $A=l^{\infty }\left( \rho
\left( \overline{w}\right) \right) $ and $B=l^{1}\left( \rho \left( 
\overline{w}\right) \right) $ $\left( \rho \in \mathcal{P}\right) $ these
functors were investigated by Janson \cite{16-Jan}. In the same situation
Janson \cite{16-Jan} also observed that (in our present notation) $%
K^{A}=Corb_{A}$ and $J^{B}=Orb_{B}.$ It is easy to see that this is also the
case for general $A$ and $B.$

\section*{1. General intermediate spaces.}

In this sec. we wish to extend Ovchinnikov's theorem \cite{23-Ovc} (as
stated in the Introduction) in the following direction. Instead of $%
l^{\infty }\left( \rho \left( \overline{w}\right) \right) $ and $l^{1}\left(
\rho \left( \overline{w}\right) \right) $ we consider quite general 
\underline{intermediate} spaces $A$ and $B$ with respect to $l^{\infty
}\left( \overline{w}\right) $ and $l^{1}\left( \overline{w}\right) $
respectively. We ask under which conditions it is true that

\begin{equation}
T:l^{\infty }\left( \overline{w}\right) \rightarrow l^{1}\left( \overline{w}%
\right) \Longrightarrow T:A\rightarrow B  \tag{1}
\end{equation}%
\emph{\ }We first derive a necessary condition, that is rather
straightforward.

\underline{Proposition 1}. Let $A$ and $B$ be intermediate spaces, with
respect to $l^{\infty }\left( \overline{w}\right) $ and $l^{1}\left( 
\overline{w}\right) $ respectively Then $\left( 1\right) $ implies that%
\begin{equation}
a=\left\{ a_{k}\right\} \in A,\dsum_{k}\frac{\left\vert b_{k}\right\vert }{%
\left\vert a_{k}\right\vert }<\infty \Longrightarrow b=\left\{ b_{n}\right\}
\in B  \tag{2}
\end{equation}%
(According to the convention made in the introduction, we interpret $\frac{%
\left\vert b_{k}\right\vert }{\left\vert a_{k}\right\vert }$ as $0$ if $%
b_{k}=0.)$

\underline{Proof:}\emph{\ }Assume $a\in A,\dsum \frac{\left\vert
b_{k}\right\vert }{\left\vert a_{k}\right\vert }<\infty .$ Let $T$ be
defined by $Tx=\left\{ \frac{b_{k}}{a_{k}}x_{k}\right\} $ if $x=\left\{
x_{k}\right\} .$($T$ is thus a multiplier transform.). Then clearly $%
T:l^{\infty }\left( \overline{w}\right) \rightarrow l^{1}\left( \overline{w}%
\right) $ and $Ta=b.$Thus $\left( 1\right) $ gives at once $b\in B.$\#$%
\medskip \medskip $

\underline{Example.}\emph{\ }If $A=l^{\infty }\left( \rho \left( \overline{w}%
\right) \right) ,B=l^{1}\left( \rho \left( \overline{w}\right) \right) $ -
Ovchinnikovs' case - then $\left( 2\right) $ is certainly is fulfilled.
Indeed if $\sup \left\vert a_{k}\right\vert <\infty ,\dsum \frac{\left\vert
b_{k}\right\vert }{\left\vert a_{k}\right\vert }<\infty $ then%
\begin{equation*}
\dsum \left\vert b_{k}\right\vert =\dsum \frac{\left\vert b_{k}\right\vert }{%
\left\vert a_{k}\right\vert }\left\vert a_{k}\right\vert \leq \dsum \frac{%
\left\vert b_{k}\right\vert }{\left\vert a_{k}\right\vert }\sup \left\vert
a_{k}\right\vert <\infty .
\end{equation*}

We now claim that the condition $\left( 2\right) $ is essentially also
necessary for $\left( 1\right) $ to hold true.

We must specify what "essentially" means.$\medskip $

To this end we make first the following general observation based on the
Aronszajn-Gagliardo theorem \cite{1-AG65}. If $A$ and $B$ are relative
interpolation spaces with respect to any Banach couples $\overline{A}$ and $%
\overline{B}$ then by enlarging $A$ and diminishing $B$ we may as well
assume that $A$ is an interpolation space with respect to $\overline{A}$ and 
$B$ an interpolation space with respect to $\overline{B}.$ Indeed we may, if
necessary, replace $A$ by $Orb_{A}\left( \overline{A}\right) \supseteq A$
and $B$ by $Corb_{B}\left( \overline{B}\right) \subseteq B;Orb_{A}\left( 
\overline{A}\right) $ and $Corb_{B}\left( \overline{B}\right) $ are then
again relative interpolation spaces with respect to $\overline{A}$ and $%
\overline{B}.$ Conversly if $Orb_{A}\left( \overline{A}\right) $ and $%
Corb_{B}\left( \overline{B}\right) $ are relative interpolation spaces with
respect to $\overline{A}$ and $\overline{B}$ so are obviously $A$ and $%
B.\medskip $

Returning to the case $\overline{A}=l^{\infty }\left( \overline{w}\right) ,%
\overline{B}=l^{1}\left( \overline{w}\right) $ we may therefore from now
assume that $A$ and $B$ are not only intermediate spaces but also \underline{%
interpolation} spaces, with respect to $l^{\infty }\left( \overline{w}%
\right) $ and $l^{1}\left( \overline{w}\right) $ respectively.

But then we have a rather explicit representation of the norm in $A$ and $B.$
In fact $A$ can be renormed in such a fashion that

\begin{equation}
\left\Vert a\right\Vert _{A}=\left\Vert \left\{ K\left( \tau _{n},a\right)
\right\} \right\Vert _{A}  \tag{3}
\end{equation}%
and $B$ in such a fashion that%
\begin{equation}
\left\Vert b\right\Vert _{B}=\inf_{b}\left\Vert \left\{ J\left( \tau
_{n},b_{n}\right) \right\} \right\Vert  \tag{4}
\end{equation}%
\emph{\ }

where the inf is extended over all representations $\widehat{b}=\left\{
b_{n}\right\} $ of $b.$We have further put $\tau _{n}=w_{n}^{0}/w_{n}^{1}.$
We also know that%
\begin{equation}
a\in l^{\infty }\left( \rho \left( \overline{w}\right) \right) \iff \sup 
\frac{K\left( \tau _{n},a\right) }{\rho \left( \tau _{n}\right) }<\infty 
\tag{5}
\end{equation}%
and that%
\begin{equation}
f\in l^{1}\left( \rho \left( \overline{w}\right) \right) \iff \dsum \frac{%
J\left( \tau _{n},b_{n}\right) }{\rho \left( \tau _{n}\right) }<\infty 
\tag{6}
\end{equation}

for some representation $\widehat{b}=\left\{ b_{n}\right\} $ on $b.$

For the proof of these facts see the Appendix.$\medskip $

It is now easy to prove

\underline{Proposition 2.} Let $A$ and $B$ be interpolation spaces with
respect to $l^{\infty }\left( \overline{w}\right) $ and $l^{1}\left( 
\overline{w}\right) .$ Then $\left( 2\right) $ implies $\left( 1\right) .$

\underline{Proof} :Let $T:l^{\infty }\left( \overline{w}\right) \rightarrow
l^{1}\left( \overline{w}\right) $ and take any $a\in A.$ We want to prove
under the hypothesis of $\left( 2\right) $ that $Ta\in B.$ We apply
Ovcinnikov's theorem \cite{23-Ovc} with $\rho \left( t\right) =K\left(
t,a\right) .$Then $a\in l^{\infty }\left( \rho \left( \overline{w}\right)
\right) $\emph{\ }by $\left( 5\right) $ so we get $Ta\in l^{1}\left( \rho
\left( \overline{w}\right) \right) .$ Therefore by $\left( 6\right) $ $Ta$
has a representation \emph{\ }$\widehat{b}$ with $\dsum \frac{J\left( \tau
_{k},b_{k}\right) }{K\left( \tau _{k},a\right) }<\infty .$ But by $\left(
3\right) $ holds $\left\{ K\left( \tau _{k},a\right) \right\} \in A$ so $%
\left( 2\right) $ gives $\left\{ J\left( \tau _{k},b_{k}\right) \right\} \in
B.$ Then by $\left( 4\right) $ $Ta\in B.$\#

\section*{2. The vector valued and the continuous cases. Retracts and
partial retracts.}

We next wish to prove a vector valued analogue of Ovchinnikov's theorem \cite%
{23-Ovc}. At the same time we will find that in the proof of Ovchinnikov's
theorem it suffices to consider geometric progressions only, that is, we can
take $\tau _{k}=\lambda ^{k}$ ("the $\lambda -adic$ case") where $\lambda $%
is a given number $>1.$(In most cases it would be sufficient to take $%
\lambda =2;$ the only raison d'\^{e}tre for allowing a general $\lambda $ is
that by choosing $\lambda $ close to $1,$ one gets the constants as close to 
$1$ as one wants to.\footnote{%
See Notes \underline{2}:$\left\langle 1\right\rangle .$}) On the other hand
in Janson's proof of this result \cite{16-Jan} an essential ingredient is
precisely the use of general ("sparse") sequences of weights. (Cf. Sec. 
\underline{$5$} ).$\medskip $

We begin with some general considerations ("abstract nonsense") on retracts
and partial retracts (cf. \cite{27-Pe}, \cite{28-Pe}), that perhaps might be
of use in other contexts too.$\medskip $

\underline{$2.1$}$.$ Let $\overline{A}$ and $\overline{X}$ be any two Banach
couples. We say that $\overline{X}$ is a \underline{partial retract} of $%
\overline{A}$ if for any $x\in \Sigma \overline{X}$ there exists linear
operators $\pi :\overline{A}\rightarrow \overline{X}$ (projection) and $%
\iota :\overline{X}\rightarrow \overline{A}$ (retraction, section of $\pi )$
such that $\pi \iota x=x.$ If we can choose $\pi $ and $\iota $
independently of $x$ we say that we have a \underline{retract} (omitting
"partial"). That is $\pi \iota =id$ and we have the commutative diagram 
\begin{subequations}
\begin{equation*}
\begin{tabular}{ccc}
$\overline{A}$ &  &  \\ 
$\downarrow \pi $ & $\nwarrow \iota $ &  \\ 
$\overline{X}$ & $\overset{id}{\rightarrow }$ & $\overline{X}$%
\end{tabular}%
\end{equation*}%
In the general case we have the commutative diagram 
\end{subequations}
\begin{equation*}
\begin{tabular}{ccc}
$\overline{A}$ &  &  \\ 
$\downarrow \pi $ & $\nwarrow \iota $ &  \\ 
$\overline{X}$ &  & $\overline{X}$ \\ 
$\uparrow \widetilde{x}$ & $\nearrow \widetilde{x}$ &  \\ 
$\overline{\mathbb{R}}$ &  & 
\end{tabular}%
\end{equation*}%
Here $\overline{\mathbb{R}}=\left( \mathbb{R},\mathbb{R}\right) $ stands for
a 1-dimensional (scalar) couple and $\widetilde{x}:c\rightarrow cx$ denotes
multiplication by $x.\medskip $

Let thus $\overline{X}$ be a partial retract of $\overline{A}$ and consider
another two Banach couples $\overline{B}$ and $\overline{Y}$ such that $%
\overline{Y}$ is partial retract of $\overline{B}.$ For $y\in \Sigma 
\overline{Y}$ denote by $\pi ^{^{\prime }}$ and $\iota ^{^{\prime }}$ the
corresponding operators.$\medskip $

The following elementary lemma is fundamental in our discussion, although it
will never be used very directly.

\underline{Lemma 1}. $y<<_{b}x\iff \iota ^{^{\prime }}y<<_{b}\iota x.$

\underline{Proof}\emph{\ : }$1):$ If $y<<_{b}x$ then $y=Tx$ for some $T:%
\overline{X}\rightarrow \overline{Y}.$ But then $\iota ^{^{\prime }}y=\iota
^{^{\prime }}T\pi \iota x=S\iota x$ where \thinspace $S=\iota ^{^{\prime
}}T\pi :\overline{A}\rightarrow \overline{B}.$ This gives $\iota ^{^{\prime
}}y<<_{b}\iota x.$\emph{\ }

$2):$ Conversely if $\iota ^{^{\prime }}y<<_{b}\iota x$ then $\iota
^{^{\prime }}y=S\iota x$ for some $S:\overline{A}\rightarrow \overline{B}.$
Now we get $y=\pi ^{^{\prime }}\iota ^{^{\prime }}y=\pi ^{^{\prime }}S\iota
x=Tx$ where $T=\pi ^{^{\prime }}S\iota :\overline{X}\rightarrow \overline{Y}%
. $ This proves $y<<_{b}x.$\# $\medskip $

\underline{Remark}\emph{\ .}The meaning of lemma $1$ is of course the
following. The fundamental question of interpolation theory is to decide for
given couples $\overline{X}$ and $\overline{Y}$ when the relation $x<<_{b}y$
is true (whatever we mean by "decide" and by "true"). So the lemma simply
says that in the case at hand this question can can be reduced to the
analogous question for the couples $\overline{A}$ and $\overline{B}$.$%
\medskip $

The following corollary illustrates the usefulness of this result.

\underline{Corollary}.\emph{\ }If $\overline{A}$ and $\overline{B}$ are
relative Calderon so are $\overline{X}$ and $\overline{Y}$.

\underline{Proof}\textrm{\ :}We have to prove $y<<_{K}x$ implies $y<<_{b}x.$%
But trivially $\iota ^{^{\prime }}y<<_{K}y$ and moreover, since $\pi \iota
x=x,x<<_{K}\iota x.$So $\iota ^{^{\prime }}y<<_{K}\iota x.$ Since $\overline{%
A}$ and $\overline{B}$ are relative Calderon this gives $\iota ^{^{\prime
}}y<<_{b}\iota x.$Lemma $1$ allows us to conclude that indeed $y<<_{b}x.$\#$%
\medskip $

We go on assuming that $\overline{X}$ is a partial retract of $\overline{A}$
and retain the notation $\pi ,\iota $ for the operators corresponding to $x$
in $\Sigma \overline{X}.$ Then we have

\underline{Lemma 2}\emph{. }Let $X$ and $A$ be relative interpolation spaces
with respect to $\overline{X}$ and $\overline{A}.$ Then $Orb_{X}\subseteq
Orb_{A}$

\underline{Proof}: Consider any Banach couple $\overline{U}$ and an element $%
u\in Orb_{X}\left( \overline{U}\right) .$We have to show that $u\in
Orb_{A}\left( \overline{U}\right) $ It is sufficient to consider the case
when $u$ is of the form $u=Tx$ with $x\in X$ and $T:\overline{X}\rightarrow 
\overline{U}.$ Then $u=T\pi \iota x=Sa$ with $S=T\pi :\overline{A}%
\rightarrow \overline{U}$ and $a=\iota x\in A$ (since $X$ and $A$ are
relative interpolation spaces with respect to $\overline{X}$ and $\overline{A%
}$ and $x\in X)$. This gives $u\in Orb_{A}\left( \overline{U}\right) .$\#$%
\medskip $

The result dual to lemma 2 is

\underline{Lemma 2$^{^{\prime }}$}. In the same hypothesis $%
Corb_{A}\subseteq Corb_{X}$

\underline{Proof}:\emph{\ }Let again $\overline{U}$ be any Banach couple and
let $u\in Corb_{A}\left( \overline{U}\right) .$We wish to show that $u\in
Corb_{X}\left( \overline{U}\right) .$ Since $u\in Corb_{A}\left( \overline{U}%
\right) $ we have $Su\in A$ for any $S:\overline{U}\rightarrow \overline{A}.$
Consider now an arbitrary operator $T:\overline{U}\rightarrow \overline{X.}$%
We have to show that $x=Tu\in X.$ Let $\iota ,\pi $ correspond to this $x$
and set $S=\iota T.$ Then clearly $S:\overline{U}\rightarrow \overline{A}$
so $Su\in A.$ But then $\pi Su\in X$, since $A$ and $X$ are relative
interpolation spaces with respect to $\overline{A}$ and $\overline{X}.$
Since $Tu=\pi Su$ this establishes our claim that $u\in Corb_{X}\left( 
\overline{U}\right) ,$\#$\medskip $

From lemma 2 and lemma $2^{^{\prime }}$we obtain the following

\underline{Corollary}\emph{. }Assume that $\overline{A}$ too is retract of $%
\overline{X}$ and that $X$ and $A$ are relative interpolation spaces with
respect to $\overline{X}$ and $\overline{A}.$ Then $%
Orb_{X}=Orb_{A},Corb_{X}=Corb_{A}.$

Thus, informally speaking, couples in the same partial retract class give
the same orbits and the same coorbits.$\medskip $

\underline{$2.2$}$.$We now return to the special case that really interest
us. We consider the Banach couple\footnote{%
See Notes \underline{2}:$\left\langle 2\right\rangle .$}%
\begin{equation*}
l^{p}\left( \overline{w}D\right) =\left( l^{p}\left( w^{0}D\right)
,l^{p}\left( w^{1}D\right) \right)
\end{equation*}%
where%
\begin{equation*}
1\leq p\leq \infty ,
\end{equation*}

\begin{equation*}
\overline{w}=\left( w^{0},w^{1}\right) =\left( \left\{ w_{n}^{0}\right\}
,\left\{ w_{n}^{1}\right\} \right)
\end{equation*}%
is a couple of weighted sequences (the index n runs through some denumerable
index set $I$), 
\begin{equation*}
D=\left\{ D_{n}\right\} \text{ is a sequence of Banach spaces }D_{n}
\end{equation*}%
\begin{equation*}
l^{p}\left( w^{i}D\right) ,\left( i=0,1\right)
\end{equation*}%
denotes the space of sequences $x=\left\{ x_{n}\right\} $ such that $%
x_{n}\in D_{n}$ for each $n\in I$ and%
\begin{equation*}
\left( \dsum_{n}\left( w_{n}^{i}\left\Vert x_{n}\right\Vert _{D_{n}}\right)
^{p}\right) ^{1/p}<\infty
\end{equation*}%
(with the usual interpretation if $p=\infty ).$

If $D_{n}=\mathbb{R}$ for each $n$ we write just $l^{p}\left( \overline{w}%
\right) .$ In particular we put%
\begin{equation*}
\overline{l}_{\lambda }^{p}=\left( l^{p}\left( \left\{ 1\right\} \right)
,l^{p}\left( \left\{ \lambda ^{-k}\right\} \right) \right) =\left(
l^{p},l^{p}\left( \lambda ^{-k}\right) \right)
\end{equation*}%
where $\lambda $ is a given number $>1;$ in this case $k$ is assumed to run
over $\mathbb{Z}$ , the set of integers. (The $\overline{}$ in the notation
is to remind us that this is couple.)$\medskip $

\underline{Proposition 1}\emph{\ . }$l^{p}\left( \overline{w}D\right) $ is a
partial retract of $\overline{l}_{\lambda }^{p}.$

\underline{Proof}: For any $k\in \mathbb{Z}$ set%
\begin{equation}
e_{k}=\left\{ n:w_{n}^{0}\leq \lambda ^{k}w_{n}^{1}<\lambda
w_{n}^{0}\right\} .  \tag{1}
\end{equation}%
\emph{\ }Clearly $e_{k}\cap e_{l}=\emptyset $ $\left( k\neq l\right) $ and $%
\cup _{k}e_{k}=I$ so we have a partition of $I$ (our index set). Let $x\in
\Sigma l^{p}\left( \overline{w}D\right) .$ We wish to exhibit the
corresponding operators $\pi :\overline{l}_{\lambda }^{p}\rightarrow
l^{p}\left( \overline{w}D\right) $ and $\iota :l^{p}\left( \overline{w}%
D\right) \rightarrow \overline{l}_{\lambda }^{p}$ with $\pi \iota x=x.$

Begin with $\iota .$ By the Hahn-Banach theorem we can find linear
functionals. $\alpha _{k}$ $\left( k\in \mathbb{Z}\right) $ on $\Sigma
l^{p}\left( \overline{w}D\right) $ such that%
\begin{equation}
\alpha _{k}\left( x\right) =\left( \dsum_{n\in e_{k}}\left(
w_{n}^{0}\left\Vert x_{n}\right\Vert _{D_{n}}\right) ^{p}\right) ^{1/p}=s_{k}
\tag{2}
\end{equation}%
and 
\begin{equation}
\left\vert \alpha _{k}\left( y\right) \right\vert \leq \left( \dsum_{n\in
e_{k}}\left( w_{n}^{0}\left\Vert y_{n}\right\Vert _{D_{n}}\right)
^{p}\right) ^{1/p}  \tag{3}
\end{equation}%
for every $y=\left\{ y_{n}\right\} \in \Sigma l^{p}\left( \overline{w}%
D\right) .$Then by the first inequality defining $e_{k}$ (see $\left(
1\right) )$ holds also%
\begin{equation}
\lambda ^{k}\left\vert \alpha _{k}\left( y\right) \right\vert \leq \left(
\dsum_{n\in e_{k}}\left( w_{n}^{1}\left\Vert y_{n}\right\Vert
_{D_{n}}\right) ^{p}\right) ^{1/p}.  \tag{4}
\end{equation}%
(Note that the expression to the right in $\left( 3\right) $ and $\left(
4\right) $ is just $\left\Vert \chi _{e_{k}}y\right\Vert _{l^{p}\left(
w^{i}D\right) }$ $\left( i=0,1\right) $ where $\chi _{e_{k}}$ denotes the
characteristic function of the set $e_{k}.).$ We set now $\iota \left(
y\right) =\left\{ \alpha _{k}\left( y\right) \right\} .$ Then $\iota
:l^{p}\left( \overline{w}D\right) \rightarrow \overline{l}_{\lambda }^{p}$.
Indeed adding up $\left( 3\right) $ and $\left( 4\right) $ we get $%
\left\Vert {}\right\Vert \leq 1.$ Also $\left( 2\right) $ gives $\iota
\left( x\right) =s=\left\{ s_{k}\right\} .$

Next let us construct $\pi .$ If $n\in e_{k}$ put $\phi _{n}=k.$ (Thus $%
n\rightarrow \phi _{n}$ is the "inverse" of the set valued function $%
k\rightarrow e_{k}.).$For $a=\left\{ a_{k}\right\} \in \Sigma \overline{l}%
_{\lambda }^{p}$ set $\pi \left( a\right) =\left\{ a_{\phi
_{n}}x_{n}/s_{\phi _{n}}\right\} .$ Then by the second inequality defining $%
e_{k}$ (see $\left( 1\right) )$ we get $\pi :\overline{l}_{\lambda
}^{p}\rightarrow l^{p}\left( \overline{w}D\right) ,\left\Vert \pi
\right\Vert \leq \lambda .$ Also clearly \thinspace $\pi \left( a\right) =s,$%
so, since $\iota \left( x\right) =s,$we finally obtain $\pi \iota x=x.$\#$%
\medskip $

In the other direction, we prove

\underline{Proposition 2\emph{.}}\emph{\ }Assume that for each $k\in \mathbb{%
Z}$ we can find $n\in e_{k}$ (where $e_{k}$ as in the proof of prop. 1, see $%
\left( 2\right) )$ such that $D_{n}\neq 0.$Then $\overline{l}_{\lambda }^{p}$
is a retract (not only a partial one) of $l^{p}\left( \overline{w}D\right) .$

\underline{Proof}: Let $n=n\left( k\right) $ be the index corresponding to $%
k $ in the hypothesis of our proposition. Pick up an element $\theta _{k}\in
D_{n\left( k\right) }$ with $\left\Vert \theta _{k}\right\Vert =1.$ We have
to construct operators $\iota :\overline{l}_{\lambda }^{p}\rightarrow
l^{p}\left( \overline{w}D\right) $ and $\pi :l^{p}\left( \overline{w}%
D\right) \rightarrow \overline{l}_{\lambda }^{p}$ such that $\pi \iota =id.$

Let $\widehat{\theta _{k}}=\left( 0,...,0,\theta _{k},0,,,\right) $ be the
sequence whose $n\left( k\right) $ th entry is $\theta _{k}$, all other
entries being zero. Then for $a=\left\{ a_{k}\right\} \in \Sigma \overline{l}%
_{\lambda }^{p}$ we set $\iota \left( a\right) =\dsum a_{k}\widehat{\theta
_{k}}/w_{n\left( k\right) }^{0}.$

The second inequality in $\left( 1\right) $ readily gives $\left\Vert \iota
\right\Vert \leq \lambda .$

Conversly select by the Hahn-Banach theorem linear functionals. on $%
D_{n\left( k\right) }$ such that $\eta _{k}\left( \theta _{k}\right)
=1,\left\Vert \eta _{k}\right\Vert =1.$Set $\pi \left( x\right) =\left\{
w_{n\left( k\right) }^{0}\eta _{k}\left( x_{n\left( k\right) }\right)
\right\} $ for $x=\left\{ x_{n}\right\} \in \Sigma l^{p}\left( \overline{w}%
D\right) .$Now the first inequality in $\left( 1\right) $ gives $\left\Vert
\pi \right\Vert \leq \lambda .$ Since clearly $\pi \iota \left( x\right) =x$
for any $x,$ the proof is complete. \#$\medskip $

\underline{2.3}\emph{\ . }After these lengthy preparations we are now ready
to prove the main result of this section.$\medskip $

We consider two couples $l^{\infty }\left( \overline{w}D\right) $ and $%
l^{1}\left( \overline{z}E\right) $ where $\overline{z}$ and $E$ have a
similar meaning as $\overline{w}$ and $D.$\emph{\ }

\underline{Theorem}.\emph{\ }$T:l^{\infty }\left( \overline{w}D\right)
\rightarrow l^{1}\left( \overline{z}E\right) \Longrightarrow T:l^{\infty
}\left( \rho \left( \overline{w}\right) D\right) \rightarrow l^{1}\left(
\rho \left( \overline{z}\right) E\right) $ for any $\rho \in \mathcal{P}$.

\underline{Proof}\emph{\ : }By \underline{$2.2$}$,$ prop. 1, $l^{\infty
}\left( \overline{w}D\right) $ is a partial retract of $\overline{l}%
_{\lambda }^{\infty }$ and $l^{1}\left( \overline{z}E\right) $ one of $%
\overline{l}_{\lambda }^{1}.$Consider any element of $l^{\infty }\left( \rho
\left( \overline{w}\right) D\right) .$We have to show that $y=Tx\in
l^{1}\left( \rho \left( \overline{z}\right) E\right) .$Let $\pi ,\iota $ and 
$\pi ^{^{\prime }},\iota ^{^{\prime }}$be the operators corresponding to $x$
and $y$ respectively. (cf. the proof of \underline{$2.1$}$,$lemma $1.)$ As
in that proof we set $S=\iota ^{^{\prime }}T\pi .$ Thus we have the
following (only "partially commutative diagram)%
\begin{equation*}
\begin{tabular}{ccc}
$l^{\infty }\left( \overline{w}D\right) $ & $\overset{T}{\rightarrow }$ & $%
l^{1}\left( \overline{z}E\right) $ \\ 
$\iota \downarrow \uparrow \pi $ &  & $\iota ^{^{\prime }}\downarrow
\uparrow \pi ^{^{\prime }}$ \\ 
$\overline{l}_{\lambda }^{\infty }$ & $\overset{S}{\rightarrow }$ & $%
\overline{l}_{\lambda }^{1}$%
\end{tabular}%
\end{equation*}%
We take for granted that $l^{\infty }\left( \frac{1}{\rho \left( \lambda
^{l}\right) }\right) $ and $l^{1}\left( \frac{1}{\rho \left( \lambda
^{k}\right) }\right) $ are relative interpolation spaces with respect to $%
\overline{l}_{\lambda }^{\infty }$ and $\overline{l}_{\lambda }^{1};$ this
is of course just the $\lambda -$adic special case of Ovchinnikov's theorem 
\cite{23-Ovc}; see introduction. It is clear that $\iota x\in l^{\infty
}\left( 1/\rho \left( \lambda ^{-k}\right) \right) .$ By the result just
quoted thus $S\iota x\in l^{1}\left( 1/\rho \left( \lambda ^{k}\right)
\right) .$It follows that $\pi ^{^{\prime }}S\iota x\in l^{1}\left( \rho
\left( \overline{z}\right) E\right) .$But $\pi ^{^{\prime }}S\iota x=\pi
^{^{\prime }}\iota ^{^{\prime }}T\pi \iota x=Tx=y.$ This completes the
proof. \#$\medskip $

\underline{Remark}\emph{\ }For a direct proof of this theorem and the
following generalization, by the method of Janson \cite{16-Jan}, see sec. 5.$%
\medskip $

\underline{2.4}.\emph{\ }We now extend the previous results $\left( 
\underline{2.2}-\underline{2.3}\right) $ to the continuous (or rather the
measurable) case.$\medskip $

Thus let $W$ be any measure space equipped with a measure $\mu .$\ We are
interested in couples%
\begin{equation*}
L^{p}\left( \overline{w}D\right) =\left( L^{p}\left( w^{0}D\right)
,L^{p}\left( w^{1}D\right) \right)
\end{equation*}%
As before $1\leq p\leq \infty .$Now $\overline{w}=\left( w^{0},w^{1}\right) $
is a couple of positive measurable functions on $W$ ("weight functions") and 
$D$ a continuous Banach bundle over $W.$ By the latter we intend a space $D$
together with a projection $\pi :D\rightarrow W$ such that each fiber $%
D_{\omega }=\pi ^{-1}\left( \omega \right) $ ($\omega \in W)$ is a Banach
space. In addition $D$ should be equipped with a "continuity structure",
i.e. there is a given vector space $\Gamma $ of sections of $D$ (called
"principal sections") such that

$\left( i\right) :$ for each $x\in \Gamma $ the map $W\rightarrow \mathbb{R}%
,\omega \rightarrow \left\Vert x\left( \omega \right) \right\Vert
_{D_{\omega }}$ is continuous,

$\left( ii\right) :$ for each $\omega \in W$ the evaluation map $\Gamma
\rightarrow D_{\omega },x\rightarrow x\left( \omega \right) $ has a dense
image in $D_{\omega }.\medskip $

Given this a continuous Banach bundle one can develop Bourbaki style
integration theory pretty much as in the special case of a trivial bundle
for which case we get the integral. In particular one can define the space
of sections $L^{p}\left( w^{i}D\right) $ $\left( i=0,1\right) $ by taking
the completion of $\Gamma $ in the norm 
\begin{equation*}
\left\Vert x\right\Vert =\left( \tint\nolimits_{W}\left( w^{i}\left( \omega
\right) \left\Vert x\left( \omega \right) \right\Vert _{D_{\omega }}\right)
^{p}d\mu \left( \omega \right) \right) ^{1/p}.
\end{equation*}

\underline{Remark}.\emph{\ }(historical) The notion of "continuous Banach
bundle" (our word) seems to be due to Godement \cite{12-God}, \cite{13-God},
the French expression being "champs continu d'espaces de Banach". Godement
develops also the corresponding integration theory.\footnote{%
See Notes \underline{2}:$\left\langle 3\right\rangle .$} (See Dixmier \cite%
{7-Dix} where ample references can be found.). If all the fibers are Hilbert
spaces and if $p=2$ one gets what is known as a "direct integral of Hilbert
spaces". (see Dixmier \cite{8-Dix}.). If $D$ is a trivial bundle, i.e. all
spaces $D_{\omega }$ coincide we take for $\Gamma $ the "constant" sections
and we are back in the classical case of vector valued integration.$\medskip 
$

Now it is easy to carry over \underline{$2.2$}$\emph{,}$ prop. 1. The
definition of $e_{k}$ now takes the form (cf. \underline{$\emph{2.2}$}, $%
\left( 1\right) )$%
\begin{equation*}
e_{k}=\left\{ \omega :w^{0}\left( \omega \right) \leq \lambda
^{k}w^{1}\left( \omega \right) <\lambda w^{0}\left( \omega \right) \right\}
\end{equation*}%
In this way we get a partition of $W$ into measurable subsets. We define
linear functionals. $\alpha _{k}$ such that $\alpha _{k}\left( x\right) =s$
with now $s_{k}=\left\Vert \chi _{e_{k}}x\right\Vert _{L^{p}\left(
w^{0}D\right) }.$The rest of the proof is pretty much the same.

At any rate \underline{we find that} $L^{p}\left( \overline{w}D\right) $ 
\underline{always is a partial retract of} $\overline{l}_{\lambda }^{p}.$ In
the same we can generalize $\underleftarrow{2.2},$ prop. $2,$ provided one
puts the proper restriction on $D$ generalizing the one of that proposition.
We will not give any details, because we are not going to use this result
anyhow.$\medskip $

Next let us consider pairs $L^{\infty }\left( \overline{w}D\right) $ and $%
L^{1}\left( \overline{z}E\right) $ where $\overline{z}$ and $E$ has a
similar meaning as $\overline{w}$ and $D;\overline{z}=\left(
z_{0},z_{1}\right) $ is thus a pair of weight functions over a measure space 
$Z,$ with measure $\nu ,$say, and $E$ is continuous Banach bundle over that
space. (Notice that since $p=\infty $ in the former case we get the same
couple $L^{\infty }\left( \overline{w}D\right) $ if we replace the measure
on $W,$ say, $\mu $ by an equivalent measure, not so in the latter case.) We
put $\rho \left( \overline{w}\right) =w^{0}/\rho \left( w^{1}/w^{0}\right) $ 
$\left( \rho \in \mathcal{P}\right) ;$ this is again a weight function on $%
W. $We define $\rho \left( \overline{z}\right) $ is a similar manner. Using
the generalization just outline (\textit{underlined}) of \underline{$2.2$}$,$
prop. 1 (the cases $p=1$ and $p=\infty $ respectively), we then get the
following analogue of \underline{$2.2$}, theorem.$\medskip $

\underline{Theorem}.\emph{\ } $T:L^{\infty }\left( \overline{w}D\right)
\rightarrow L^{1}\left( \overline{z}E\right) \Longrightarrow T:L^{\infty
}\left( \rho \left( \overline{w}\right) D\right) \rightarrow L^{1}\left(
\rho \left( \overline{z}\right) E\right) $ for any $\rho \in \mathcal{P}$ .$%
\medskip $

\underline{$2.5$}$.$ Let us indicate an application of the previous theorem
. Consider the \underline{Besov Couple}%
\begin{equation*}
B_{p}^{\overline{s}r}\left( \mathbb{R}^{n}\right) =\left(
B_{p}^{s_{0}r}\left( \mathbb{R}^{n}\right) ,B_{p}^{s_{1}r}\left( \mathbb{R}%
^{n}\right) \right)
\end{equation*}%
with $\overline{s}=\left( s_{0},s_{1}\right) $ and $r,p\in \left[ 1,\infty %
\right] .$ It is well-known (see e.g. \cite[chap. 6]{3-BeL0} or \cite[chap. 5%
]{29-Pe}) that $B_{p}^{\overline{s}r}\left( \mathbb{R}^{n}\right) $ is a
retract of $l^{p}\left( 2^{k\overline{s}}L_{p}\right) $ where $l^{p}\left(
2^{k\overline{s}}\right) $ denotes the couple of weighted sequences $%
l^{p}\left( \left\{ 2^{ks_{0}}\right\} ,\left\{ 2^{ks_{1}}\right\} \right) $
and $L_{p}=L_{p}\left[ 0,1\right] .$Using this fact \underline{$2.4$}$,$
theorem thus yields as a corollary e.g. the following result.$\medskip $

\underline{Theorem}.\emph{\ }$T:B_{p}^{\overline{s}\infty }\left( \mathbb{R}%
^{n}\right) \rightarrow B_{q}^{\overline{t}1}\left( \mathbb{R}^{m}\right)
\Longrightarrow T:B_{p}^{s\infty }\left( \mathbb{R}^{n}\right) \rightarrow
B_{q}^{t1}\left( \mathbb{R}^{m}\right) $ provided $s=\left( 1-\theta \right)
s_{0}+\theta s_{1},t=\left( 1-\theta \right) t_{0}+\theta t_{1},\theta \in
\left( 0,1\right) .$ Here $\overline{s}=\left( s_{0},s_{1}\right) ,\overline{%
t}=\left( t_{0},t_{1}\right) ,p,q\in \left[ 1,\infty \right] $ are arbitrary.%
$\medskip $

\underline{Example}.\emph{\ }Taking $m=n$ this result is formally applicable
to the convolution operator $Tf=a\ast f$ where $a\in B_{r}^{\gamma r}$ where 
$\frac{1}{r}=\frac{1}{p}+\frac{1}{q}-1$ and $t_{i}=s_{i}+\gamma $ $\left(
i=0,1\right) .$The conclusion in this case follows of course also directly,
without any recourse to such a sophisticated device as Ovchinnikov's theorem 
\cite{23-Ovc}. \underline{This is perhaps symptomatic}\textit{. }At any rate
we - this is of course no proof - do not know of any non-trivial application
of this result to the usual operators of Analysis.

\section*{3. Condition $\left( O\right) .$Type $\left( l\right) $ and type $%
\left( u\right) .$}

Now we wish to pass to a more abstract situation.$\medskip $

\underline{3.1.} Let us first return to the primitive situation of sec. 
\underline{1}. (We have not been able to extend the results of that sec. to
the vector valued case of sec. \underline{2}.). In view of the results of
sec. \underline{2}. it is no great case of generality to restrict oneself to
the $\lambda -$adic case, i.e. we consider (in the notation of that sec.)
the couples $\overline{l}_{\lambda }^{\infty }$ and $\overline{l}_{\lambda
}^{1}.$ If $A$ and $B$ are relative interpolation spaces with respect to $%
\overline{l}_{\lambda }^{\infty }$ and $\overline{l}_{\lambda }^{1}$ then we
must have the inclusion $G^{A}\subseteq H^{B}.$ If in addition $A$ is an
interpolation space with respect to $\overline{l}_{\lambda }^{\infty }$ and $%
B$ one with respect to $\overline{l}_{\lambda }^{1}$ then by Janson \cite%
{16-Jan} we get $G^{A}\left( \overline{l}_{\lambda }^{\infty }\right)
=K^{A}\left( \overline{l}_{\lambda }^{\infty }\right) =A,H^{B}\left( 
\overline{l}_{\lambda }^{1}\right) =J^{B}\left( \overline{l}_{\lambda
}^{1}\right) =B.$ (For the definition of the functors $%
G^{A},K^{A},H^{B},J^{B}$ sec. \underline{$\frac{1}{2}$}. ). This leads to
the following question: Given any two Banach couples $\overline{X}$ and $%
\overline{Y}$ under which conditions can we assert that $T:\overline{X}%
\rightarrow \overline{Y}$ implies $T:K^{A}\left( \overline{X}\right)
\rightarrow J^{B}\left( \overline{Y}\right) $ for any $A$ and $B,$i.e. that $%
K^{A}\left( \overline{X}\right) $ and $J^{B}\left( \overline{Y}\right) $ are
relative interpolation spaces with respect to $\overline{X}$ and $\overline{Y%
}.$ (That $T:\overline{X}\rightarrow \overline{Y}$ implies $T:G^{A}\left( 
\overline{X}\right) \rightarrow H^{B}\left( \overline{Y}\right) $ is by what
we just said obvious.). We then say that $\overline{X}$ \underline{and} $%
\overline{Y}$ \underline{satisfy condition } $\left( O\right) .$ $\medskip $

We have thus already met (see sec. \underline{2}) several instances of
couples $\overline{X}$ and $\overline{Y}$ satisfying condition $\left(
O\right) ,\overline{X}=\overline{l}^{\infty },\overline{Y}=\overline{l}^{1}$
being the primitive case, more generally $\overline{X}=l^{\infty }\left( 
\overline{w},D\right) ,\overline{Y}=l^{1}\left( \overline{z},E\right) $ or $%
\overline{X}=L^{\infty }\left( \overline{w},D\right) ,\overline{Y}%
=L^{1}\left( \overline{z},E\right) ,$ likewise $\overline{X}=B_{p}^{%
\overline{s},\infty }\left( \mathbb{R}^{n}\right) ,\overline{Y}=B_{q}^{%
\overline{t},1}\left( \mathbb{R}^{m}\right) .$ As is easy to convince
oneself various combinations of these couples will also do.$\medskip $

Let us therefore coin the following additional terminology. We say that $%
\overline{X}$ \underline{is of type } $\left( l\right) $ if $K^{A}\left( 
\overline{X}\right) =G^{A}\left( \overline{X}\right) $ for any interpolation
space $A$ with respect to $\overline{l}_{\lambda }^{\infty }$ and that $%
\overline{Y}$ \underline{is of type } $\left( u\right) $ if $H^{B}\left( 
\overline{Y}\right) =J^{B}\left( \overline{Y}\right) $ for any interpolation
space $B$ with respect to $\overline{l}_{\lambda }^{1}.$ E.g., $\overline{l}%
_{\lambda }^{r},l^{r}\left( \overline{w},D\right) ,L^{r}\left( \overline{w}%
,D\right) ,B_{p}^{\overline{s},q}\left( \mathbb{R}^{n}\right) $ are of type $%
\left( l\right) $ if $r=\infty $ and of type $\left( u\right) $ if $r=1.$%
(The letters l and u are of course chosen to honour Ovchinnikov.)$\medskip $

The following result is almost trivial. (It has of course been implicit the
above discussion.)$\medskip $

\underline{Proposition. } Let $\overline{X}$ be of type $\left( l\right) $
and $\overline{Y}$ of type $\left( u\right) .$Then $\overline{X}$ and $%
\overline{Y}$ satisfy condition $\left( O\right) .$

\underline{Proof }: Let $T:\overline{X}\rightarrow \overline{Y}.$Then $%
T:G^{A}\left( \overline{X}\right) \rightarrow H^{B}\left( \overline{Y}%
\right) $ since $G^{A}\subseteq H^{B}$ whenever $A$ and $B$ are relative
interpolation spaces with respect to $\overline{l}_{\lambda }^{\infty }$ and 
$\overline{l}_{\lambda }^{1}.$ But $K^{A}\left( \overline{X}\right)
=G^{A}\left( \overline{X}\right) ,H^{B}\left( \overline{Y}\right)
=J^{B}\left( \overline{Y}\right) $ by our assumption on $\overline{X}$ and $%
\overline{Y}.$ Se we get indeed $T:K^{A}\left( \overline{X}\right)
\rightarrow J^{B}\left( \overline{Y}\right) .$\#$\medskip $

\underline{Remark}. Indeed we know of no single non-trivial instance of two
couples $\overline{X}$ and $\overline{Y}$ satisfying condition $\left(
O\right) $ where $\overline{X}$ is not of type $\left( l\right) $ and $%
\overline{Y}$ not of type $\left( u\right) .\medskip $

\underline{3.2}. In proving that two given Banach couples satisfy condition $%
\left( O\right) $ it is not necessary to consider the most general
interpolation spaces with respect to $\overline{l}_{\lambda }^{\infty }$ and 
$\overline{l}_{\lambda }^{1}.$ Indeed the following holds true.$\medskip $

\underline{Proposition}. Let $\overline{X}$ and $\overline{Y}$ be any two
Banach couples . The following conditions are equivalent.

$\left( i\right) :$ $\overline{X}$ and $\overline{Y}$ satisfy condition $%
\left( O\right) .$

$\left( ii\right) :T:\overline{X}\rightarrow \overline{Y}\Longrightarrow T:%
\overline{X}_{\rho \infty :K}\rightarrow \overline{Y}_{\rho 1:J}$ for any $%
\rho \in \mathcal{P}$.

$\left( iii\right) :$ For any $y\in \Sigma _{0}\overline{Y}$ and $x\in
\Sigma \overline{X}$ holds $y<<_{b}x\Longrightarrow y<<_{J/K}x.$

\underline{Proof}. The implications $\left( i\right) \Longrightarrow \left(
ii\right) \Longrightarrow \left( iii\right) $ are trivial. There remains the
implication $\left( iii\right) \Longrightarrow \left( i\right) ,$ which is
proved just by adopting the proof of \underline{1}, prop. 2. Indeed assume
that $x\in K^{A}\left( \overline{X}\right) $ i.e. $\left\{ K\left( \lambda
^{k},x\right) \right\} \in A$ and let there be given a bounded linear
operator $T:\overline{X}\rightarrow \overline{Y}.$ We must show that $%
y=Tx\in J^{B}\left( \overline{Y}\right) .$ (Here $A$ and $B$ are relative
interpolation spaces with respect to $\overline{l}_{\lambda }^{\infty }$ and 
$\overline{l}_{\lambda }^{1}.).$ Since $y<<_{b}x$ the condition $\left(
iii\right) $ gives us a representation $\widehat{y}$ of $y$ such that $%
\Sigma J\left( \lambda ^{k},y_{k}\right) /K\left( \lambda ^{k},x\right)
<\infty .$ Now, by our assumption on $A$ and $B$ and by $\underline{1},$
prop. \underline{1}, condition 1, $\left( 2\right) $ must be true. Therefore
we conclude $\left\{ J\left( \lambda ^{k},y_{k}\right) \right\} \in B.$
Hence $y\in J^{B}\left( \overline{Y}\right) .$\#$\medskip $

\underline{Corollary}. If $\overline{X}$ and $\overline{Y}$ satisfy
condition $\left( O\right) $ then $\overline{Y}$ must be regular $\left(
\Sigma _{0}\overline{Y}=\Sigma \overline{Y}\right) .\medskip $

\underline{3.3}. Next we give useful criteria for a given couple to be type $%
\left( l\right) $ or of type $\left( u\right) .$

We begin with the following. $\medskip $

\underline{Proposition 1}. Let $\overline{X}$ and $\overline{X}^{\left(
1\right) }$ be any two Banach couples. Assume that for each $x\in \Sigma 
\overline{X}$ there exists an element $x^{\left( 1\right) }\in \Sigma 
\overline{X}^{\left( 1\right) }$ and a linear operator $\pi :\overline{X}%
^{\left( 1\right) }\rightarrow \overline{X}$ such that $x=\pi x^{\left(
1\right) }$ and $x^{\left( 1\right) }<<_{K}x.$ If $\overline{X}^{\left(
1\right) }$ is of type $\left( l\right) $ so is $\overline{X}.$

\underline{Proof}. It is convenient to picture the situation diagrammatically%
\begin{equation*}
\begin{tabular}{cc}
& $\overline{X}$ \\ 
$\pi $ & $\uparrow $ \\ 
& $\overline{X}^{\left( 1\right) }$%
\end{tabular}%
,%
\begin{tabular}{cc}
$x$ &  \\ 
$\uparrow $ & $<<_{K}$ \\ 
$x^{\left( 1\right) }$ & 
\end{tabular}%
\end{equation*}%
Assuming that $x\in K^{A}\left( \overline{X}\right) ,A$ an interpolation
space with respect to $\overline{l}_{\lambda }^{\infty },$ we have to show
that $x\in G^{A}\left( \overline{X}\right) .$ Now $\left\{ K\left( \lambda
^{j},x\right) \right\} \in A,$ since $x\in K^{A}\left( \overline{X}\right) ,$
so $x^{\left( 1\right) }<<_{K}x$ gives $\left\{ K\left( \lambda
^{j},x^{\left( 1\right) }\right) \right\} \in A.$ Thus $x^{\left( 1\right) }$
is a linear combinations of elements of the form $Ta$ where $a\in A$ and $T:%
\overline{l}_{\lambda }^{\infty }\rightarrow \overline{X}^{\left( 1\right)
}. $ Then $x$ is linear combination of elements $T^{^{\prime }}a$ where $%
a\in A$ and $T^{^{\prime }}=\pi T:\overline{l}_{\lambda }^{\infty
}\rightarrow \overline{X}.$ This proves $x\in G^{A}\left( \overline{X}%
\right) .$ Thus $K^{A}\left( \overline{X}\right) =G^{A}\left( \overline{X}%
\right) $ (for any $A)$ and $\overline{X}$ is effectively of type $\left(
l\right) .$\#$\medskip $

Next we prove the dual result.

\underline{Proposition 2}. Let $\overline{Y}$ and $\overline{Y}^{\left(
1\right) }$ be any two Banach couples. Assume that for each $y\in \Sigma _{0}%
\overline{Y}$ and a linear operator $\iota :\overline{Y}\rightarrow 
\overline{Y}^{\left( 1\right) }$ such that $y^{\left( 1\right) }=\iota y$
and $y<<_{J}y^{\left( 1\right) }.$ If $\overline{Y}^{\left( 1\right) }$ is
of type $\left( u\right) $ then so is $\overline{Y}$.

\underline{Proof}. The proof is parallel to the one of prop. 1, is possibly
still simpler. The relevant diagram is now%
\begin{equation*}
\begin{tabular}{cc}
& $\overline{Y}$ \\ 
$\iota $ & $\downarrow $ \\ 
& $\overline{Y}^{\left( 1\right) }$%
\end{tabular}%
,%
\begin{tabular}{cc}
$y$ &  \\ 
$\downarrow $ & $<<_{J}$ \\ 
$y^{\left( 1\right) }$ & 
\end{tabular}%
\end{equation*}%
Let $y\in H^{B}\left( \overline{Y}\right) ,$ $B$ an interpolation space with
respect to $\overline{l}_{\lambda }^{1}.$ We must verify that $y\in
J^{B}\left( \overline{Y}\right) .$ Since $y\in H^{B}\left( \overline{Y}%
\right) $ for any $T:\overline{Y}\rightarrow \overline{l}_{\lambda }^{1}$ we
have $Ty\in B.$But then in particular $T^{1}\iota y\in B$ for any $T^{1}:%
\overline{Y}^{\left( 1\right) }\rightarrow \overline{l}_{\lambda }^{1},$
that is $T^{1}y^{\left( 1\right) }\in B.$ Thus, by the same token as above,
we conclude $y^{\left( 1\right) }\in H^{B}\left( \overline{Y}^{\left(
1\right) }\right) $ or, since $H^{B}\left( \overline{Y}^{\left( 1\right)
}\right) =J^{B}\left( \overline{Y}^{\left( 1\right) }\right) $ by
hypothesis, $y^{\left( 1\right) }\in J^{B}\left( \overline{Y}^{\left(
1\right) }\right) .$ Since $y<<_{J}y^{\left( 1\right) }$ this gives $y\in
J^{B}\left( \overline{Y}\right) .$ We have proved that $H^{B}\left( 
\overline{Y}\right) =J^{B}\left( \overline{Y}\right) $ (for any $B)$ and $%
\overline{Y}$ is of type $\left( u\right) .$\#$\medskip $

\underline{Remark}. In the somewhat obscure terminology of \cite{27-Pe} the
hypothesis of prop. 2 says that $\overline{Y}$ is a lb-pseudoretract of $%
\overline{Y}^{\left( 1\right) }.$ Similarly the condition of prop. 1 means
that $\overline{X}$ is a lb-pseudoretract of $\overline{X}^{\left( 1\right)
}.$ In sec. \underline{2} of the present paper we further used "partial
retract" in place for lb-pseudoretract. (b for "bounded", l for "linear".)$%
\medskip $

In praxis the relation $<<_{J}$ is difficult to verify (except in trivial
cases). Luckily we can in the special case of interest for us (cf. infra 
\underline{3.4}) substitute it for $<<_{K}.$Namely there holds the following.

\underline{Corollary}. Assume that $\overline{Y}$ satisfies the "strong form
of the fundamental lemma" (see \cite{6-CwPe}, p. 33, condition $\left(
3\right) $). Then we can make the same conclusion as in prop. 2 also if we
only assume $y<<_{K}y^{\left( 1\right) }$ (not necessarily $y<<_{J}y^{\left(
1\right) }).$

$\medskip $

Before proving this result we make some clarifications.

\underline{Remark}. Let $\overline{A}$ be any Banach couple. The
"fundamental lemma" (see e.g. \cite{3-BeL0}, p. 45) says that any $a\in
\Sigma _{0}\overline{A}$ has a representation $\widehat{a}=\left\{
a_{k}\right\} _{k\in \mathbb{Z}}$ such that $J\left( \lambda
^{k},a_{k}\right) \leq CK\left( \lambda ^{k},a\right) $ for all $k\in 
\mathbb{Z}$ where $C$ is a constant depending in $\lambda $ only. The
"strong form of the fundamental lemma" (see loc. cit.) says that one can
choose $\widehat{a}$ such that $\Sigma _{j}\min \left( 1,\frac{\lambda ^{k}}{%
\lambda ^{j}}\right) J\left( \lambda ^{j},a_{j}\right) \leq C^{^{\prime
}}K\left( \lambda ^{k},a\right) $ (or in brief $\Omega \left\{ J\left(
\lambda ^{j},a_{j}\right) \right\} \leq C^{^{\prime }}K\left( \lambda
^{k},a\right) ;\Omega $ is the Calderon operator, see the appendix) for all $%
k\in \mathbb{Z}$ where $C^{^{\prime }}$ now depends on $\overline{A}$ too.
What couples $\overline{A}$ admit this strong form of the fundamental lemma
is not quite clear yet (cf. \cite{3-BeL0}).$\medskip $

After this digression we proceed to the

\underline{Proof} (of the corollary of prop. 2): Assume that $%
y<<_{K}y^{\left( 1\right) },y^{\left( 1\right) }\in J^{B}\left( \overline{Y}%
^{\left( 1\right) }\right) .$We have to prove that $y\in J^{B}\left( 
\overline{Y}\right) .$ Let the norm in $B$ be given by%
\begin{equation*}
\left\Vert b\right\Vert _{B}\approx \Phi \left( \left\{ K\left( \lambda
^{k},b\right) \right\} \right) \approx \inf_{\widehat{b}}\Phi \left( \Omega
\left( \left\{ J\left( \lambda ^{k},b_{k}\right) \right\} \right) \right) 
\end{equation*}%
where $\Phi $ is a sequence norm and (in the last expression) $\Omega $
stands for the Calderon transformation, the inf being taken over all
representations $\widehat{b}$ of $b$ (see the appendix). Then we must have $%
\Phi \left( \Omega \left( \left\{ J\left( \lambda ^{k},y_{k}^{\left(
1\right) }\right) \right\} \right) \right) <\infty $ for some representation 
$\widehat{y}^{\left( 1\right) }=\left\{ y_{k}^{\left( 1\right) }\right\} $
of $y^{\left( 1\right) }.$ Since always $K\left( \lambda ^{k},y^{\left(
1\right) }\right) \leq \Omega \left\{ J\left( \lambda ^{k},y_{k}^{\left(
1\right) }\right) \right\} $ (see e.g. \cite{3-BeL0}, p.44) this gives $\Phi
\left( \left\{ K\left( \lambda ^{k},y^{\left( 1\right) }\right) \right\}
\right) <\infty $ and thus a fortiori - remember that $y<<_{K}y^{\left(
1\right) }$- $\Phi \left( \left\{ K\left( \lambda ^{k},y\right) \right\}
\right) <\infty .$ But if $\widehat{y}$ is the representation of $y$
provided by the strong form of the fundamental lemma then also $\Phi \left(
\Omega \left( \left\{ J\left( \lambda ^{k},y_{k}\right) \right\} \right)
\right) <\infty .$ This proves that $y\in J^{B}\left( \overline{Y}\right) .$%
\#$\medskip $

\underline{3.4}. As an application of the proceeding considerations we now
show that Marcinkiewicz and Lorentz couples are of type $\left( l\right) $
and of type $\left( u\right) $ respectively. $\medskip $

Let us first recall the definition of Marcinkiewicz $M\left( \psi \right) $
and Lorentz spaces $\Lambda \left( \phi \right) .$ (A good introduction to
Marcinkiewicz and Lorentz spaces can be found in chap. 2 of the book by
Krein-Petunin-Semenov \cite{17-KrPeSe}.)$\medskip $

In what follows $\psi $ and $\phi $ denote positive concave functions on $%
\left( 0,\infty \right) $ $\left( \psi ,\phi \in \mathcal{K}\right) $ such
that%
\begin{equation}
\tint\nolimits_{0}^{t}\frac{du}{\psi \left( u\right) }\leq C_{1}\frac{t}{%
\psi \left( t\right) }  \tag{1}
\end{equation}%
and%
\begin{equation}
\tint\nolimits_{0}^{t}\phi \left( u\right) \frac{du}{u}\leq C_{2}\phi \left(
t\right)  \tag{2}
\end{equation}%
respectively, with suitable constants $C_{1}$ and $C_{2}$.$\medskip $

We consider a fixed otherwise unspecified measure space. (Below we shall
specialize to the case of the interval $\left( 0,\infty \right) $ equipped
with the measure $dt.$). We say that a measurable function $f$ on this space
belongs to $M\left( \psi \right) $ if 
\begin{equation}
\sup_{t}\frac{\psi \left( t\right) }{t}\tint\nolimits_{0}^{t}f^{\ast }\left(
u\right) du<\infty  \tag{3}
\end{equation}%
or equivalently $\sup_{t}\psi \left( t\right) f^{\ast }\left( t\right)
<\infty $ and to $\Lambda \left( \phi \right) $ if 
\begin{equation}
\tint\nolimits_{0}^{\infty }f^{\ast }\left( t\right) d\phi \left( t\right)
<\infty  \tag{4}
\end{equation}%
or equivalently $\tint\nolimits_{0}^{\infty }f^{\ast }\left( t\right) \phi
\left( t\right) \frac{dt}{t}<\infty .$ Here $f^{\ast }$ stands for the
non-increasing rearrangement of $f.$ The equivalence of the two conditions
in $\left( 3\right) $ or $\left( 4\right) $ results from $\left( 1\right) $
and $\left( 2\right) $ respectively.\footnote{%
See Notes \underline{3}:$\left\langle 1\right\rangle $}$\medskip $

One can show that (see \cite{17-KrPeSe})%
\begin{equation}
M\left( \psi _{0}\right) +M\left( \psi _{1}\right) =M\left( \min \left( \psi
_{0},\psi _{1}\right) \right)  \tag{5}
\end{equation}%
and that%
\begin{equation}
\Lambda \left( \phi _{0}\right) +\Lambda \left( \phi _{1}\right) =\Lambda
\left( \min \left( \phi _{0},\phi _{1}\right) \right)  \tag{6}
\end{equation}%
up to equivalence of norm. Replacing here $\psi _{1}$ and $\phi _{1}$ by $%
t\psi _{1}$ and $t\phi _{1}$ we get estimates for the K-functional in the $%
\underline{\text{Marcinkiewicz couple}}$ $M\left( \overline{\psi }\right)
=\left( M\left( \psi _{0}\right) ,M\left( \psi _{1}\right) \right) $ and the 
\underline{Lorentz couple} $\Lambda \left( \overline{\phi }\right) =\left(
\Lambda \left( \phi _{0}\right) ,\Lambda \left( \phi _{1}\right) \right) .$
(One has similar estimates for the J-functional (see \cite{17-KrPeSe}) but
we do not need them here.).$\medskip $

We can now announce the following result

\underline{Theorem}. $\left( i\right) $ The Marcinkiewicz couple $M\left( 
\overline{\psi }\right) =\left( M\left( \psi _{0}\right) ,M\left( \psi
_{1}\right) \right) $, where each $\psi _{i}$ $\left( i=0,1\right) $ is
subject to $\left( 1\right) $ is of type $\left( l\right) .$ $\left(
ii\right) $ The Lorentz couple $\Lambda \left( \overline{\phi }\right)
=\left( \Lambda \left( \phi _{0}\right) ,\Lambda \left( \phi _{1}\right)
\right) ,$where each $\phi _{i}$ $\left( i=0,1\right) $ is subject to $%
\left( 2\right) ,$ is of type $\left( u\right) .$

\underline{Proof}. For simplicity we give the proof only for the case when
the underlying measure space is the interval $\left( 0,\infty \right) $
equipped with the Lebesque measure $dt.$

\underline{Ad $\left( i\right) $} . We wish to apply \underline{3.3} prop 1
with $\overline{X}=\left( M\left( \psi _{0}\right) ,M\left( \psi _{1}\right)
\right) $ and $\overline{X}^{\left( 1\right) }=\left( L^{\infty }\left( \psi
_{0}\right) ,L^{\infty }\left( \psi _{1}\right) \right) .$ To a given $f$ in 
$\Sigma X$ we have to produce a suitable "projection" $\pi :\overline{X}%
^{\left( 1\right) }\rightarrow \overline{X}.$ With no loss of generality we
may assume that $f=f^{\ast }.$ We then let $\pi $ be simply the identity
map. Its continuity is established as follows. Let $g$ be an element of $%
L^{\infty }\left( \psi \right) $ where $\psi $ is any function (in $\mathcal{%
K}$) subject to $\left( 1\right) .$Then $\left\vert g\left( t\right)
\right\vert \leq \frac{C}{\psi \left( t\right) }$ a.e. on $\left( 0,\infty
\right) $ for some constant $C.$ Now $\frac{1}{\psi }$ is a non-increasing
function. Thus we get $g^{\ast }\left( t\right) \leq \frac{C}{\psi \left(
t\right) }.$ This shows (cf. $\left( 1\right) )$ that $g\in M\left( \psi
\right) .$ Thus the identity map from $L^{\infty }\left( \psi \right) $ into 
$M\left( \psi \right) $ is continuous. In particular taking $\psi =\psi _{i}$
$\left( i=0,1\right) $ this substantiates our previous claim. Since our $f$
has the same K-functional in both couples, up to equivalence%
\begin{equation*}
K\left( t,f\right) \approx \sup \min \left( \psi _{0}\left( u\right) ,t\psi
_{1}\left( u\right) \right) f^{\ast }\left( u\right) 
\end{equation*}%
- here we use $\left( 5\right) $- all the conditions for the application of 
\underline{3.3} prop. 1 are met. We conclude that, since $\overline{X}%
^{\left( 1\right) }=L^{\infty }\left( \overline{\psi }\right) $ is of type $%
\left( l\right) $ (see \underline{3.2} ), so is $\overline{X}=M\left( 
\overline{\psi }\right) .$

\underline{Ad $\left( ii\right) $}. This time we wish to apply \underline{3.2%
} cor. of prop. 2 taking $\overline{Y}=\Lambda \left( \overline{\phi }%
\right) ,\overline{Y}^{\left( 1\right) }=L^{1}\left( \frac{\overline{\phi }}{%
t}\right) $ (where $\frac{\overline{\phi }}{t}$ is the couple of weight
functions $\left( \frac{\phi _{0}\left( t\right) }{t},\frac{\phi _{1}\left(
t\right) }{t}\right) ).$ Again we know that this $\overline{Y}^{\left(
1\right) }$ is of type $\left( u\right) $ (see \underline{3.3}). That $%
\overline{Y}$ admits the strong form of the fundamental lemma is indicated
in \cite{6-CwPe}, p. 33. To a given $f$ in $\Sigma \overline{Y}$ we have to
produce a suitable "retraction" $\iota :\overline{Y}\rightarrow \overline{Y}%
^{\left( 1\right) }$ and again we fix attention to the case $f=f^{\ast }$ in
which case we can take $\iota $ to be simply the identity map. The proof of
its continuity runs as follows. Consider quite generally $g$ in $\Lambda
\left( \phi \right) $ when $\phi $ is any function (in $\mathcal{K)}$
subject to $\left( 2\right) .$Then by $\left( 2\right) $ $%
\tint\nolimits_{0}^{\infty }g^{\ast }\left( t\right) \frac{\phi \left(
t\right) }{t}<\infty .$But $\frac{\phi \left( t\right) }{t}$ is a
non-increasing function. Thus by a basic property of rearrangements (see
e.g. \cite{3-BeL0}, p. 5) it follows that $\tint\nolimits_{0}^{\infty
}\left\vert g\left( t\right) \right\vert \frac{\phi \left( t\right) }{t}%
dt<\infty .$This proves $g\in L^{1}\left( \frac{\phi \left( t\right) }{t}%
\right) .$ In particular taking $\phi =\phi _{i}$ $\left( i=0,1\right) $ we
get the desired continuity of our map $\iota .$ For the K-functional of $f$
we have in both couples the estimate (by $\left( 6\right) $ \underline{and } 
$\left( 2\right) )$%
\begin{equation*}
K\left( t,f\right) \approx \tint\nolimits_{0}^{\infty }f^{\ast }\left(
u\right) \min \left( \phi _{0}\left( u\right) ,t\phi _{1}\left( u\right)
\right) \frac{du}{u}
\end{equation*}%
Thus \underline{3.2} cor. of prop. 2 is effectively applicable and we
conclude that our $\Lambda \left( \overline{\phi }\right) $ indeed is of
type $\left( u\right) .$\#$\medskip $

The following result is now immediate (formal application of \underline{3.1}
prop.)

\underline{Corollary}. Consider any Marcinkiewicz couple $M\left( \overline{%
\psi }\right) $ and any Lorentz couple $\Lambda \left( \overline{\phi }%
\right) $ (possibly over \underline{different} measure spaces). Then $%
M\left( \overline{\psi }\right) $ and $\Lambda \left( \overline{\phi }%
\right) $ satisfy condition $\left( O\right) .$ In particular thus if $%
T:M\left( \overline{\psi }\right) \rightarrow \Lambda \left( \overline{\phi }%
\right) $ then $T:M\left( \rho \left( \overline{\psi }\right) \right)
\rightarrow \Lambda \left( \rho \left( \overline{\phi }\right) \right) $ for
any $\rho \in \mathcal{P}$.$\medskip $

\underline{Remark}. As in the case of Besov couples (sec. \underline{2.5},
remark) we know of no interesting operators mapping a space $M\left( \psi
\right) $ into space $\Lambda \left( \phi \right) .$

\section*{4. Nuclearity.}

\underline{4.1}. The following result, which is fundamental for our
discussion, appears in the literature in many guises It is often associated
with the names of Bergh and Cwikel. See the following remark (historical).

\underline{Lemma}. If $y<_{n}x$ then $y<_{J/K\left( \lambda \right) }\lambda
x.$ Conversely if $y<_{J/K\left( \lambda \right) }x$ then $y<_{n}x$

\underline{Proof}: Assume $y<_{n}x.$ Then by definition $y=Tx$ for some $T:%
\overline{X}\overset{n}{\rightarrow }\overline{Y}$ with $\left\Vert
T\right\Vert _{n}<1+\epsilon ,\epsilon >0,$ i.e. $Tx=\dsum_{n}l_{n}\left(
a\right) b_{n}$ with $\dsum_{n}\max_{i=0,1}\left\Vert l_{n}\right\Vert
_{X_{i}^{^{\prime }}}\left\Vert b_{n}\right\Vert _{Y_{i}}<1+\epsilon .$ Let $%
e_{k}$ $\left( k\in \mathbb{Z}\right) $ be the subset of the our
(denumerable) index set $I,$ say defined by the condition%
\begin{equation*}
\left\Vert b_{n}\right\Vert _{Y_{0}}\leq \lambda ^{k}\left\Vert
b_{n}\right\Vert _{Y_{1}}<\lambda \left\Vert b_{n}\right\Vert _{Y_{0}}.
\end{equation*}%
Then%
\begin{equation}
J\left( \lambda ^{-k},l_{n}\right) J\left( \lambda ^{k},b_{n}\right) \leq
\lambda \max_{i=0,1}\left\Vert l_{n}\right\Vert _{X_{i}^{^{\prime
}}}\left\Vert b_{n}\right\Vert _{Y_{i}}.  \tag{1}
\end{equation}%
With no loss of generality we may assume $b_{n}\neq 0$ for all $n.$So each $%
n $ belongs to exactly one $e_{k}$ and we have a partition of $I:I=\cup
_{k\in \mathbb{Z}}e_{k},e_{k}\cap e_{l}=\emptyset $ if $k\neq l.$ If we set $%
y_{k}=\dsum_{n\in e_{k}}l_{n}\left( x\right) b_{n}$ (the series is clearly
summable in $\Sigma \overline{Y})$ we therefore have $y=\dsum_{k\in \mathbb{Z%
}}y_{k}$, that is $\widehat{y}=\left\{ y_{k}\right\} _{k\in \mathbb{Z}}$ is
a representation of $y.$Furthermore by $\left( 1\right) $%
\begin{eqnarray*}
J\left( \lambda ^{k},y_{k}\right) &\leq &\dsum_{n\in e_{k}}J\left( \lambda
^{-k},l_{n}\right) K\left( \lambda ^{k},x\right) J\left( \lambda
^{k},b_{n}\right) \\
&\leq &\lambda \dsum_{n\in e_{k}}\max_{i=0,1}\left\Vert l_{n}\right\Vert
_{X_{i}^{^{\prime }}}\left\Vert b_{n}\right\Vert _{Y_{i}}K\left( \lambda
^{k},x\right)
\end{eqnarray*}%
which yields%
\begin{equation*}
\dsum_{k}\frac{J\left( \lambda ^{k},y_{k}\right) }{K\left( \lambda
^{k},x\right) }\leq \lambda \dsum_{n}\max_{i=0,1}\left( \left\Vert
l_{n}\right\Vert _{X_{i}^{^{\prime }}}\left\Vert b_{n}\right\Vert
_{Y_{i}}\right) <\lambda \left( 1+\epsilon \right) .
\end{equation*}%
Since $\epsilon >0$ is arbitrary, this shows that $y\leq _{J/K\left( \lambda
\right) }\lambda x.$

Conversly let $y<_{J/K\left( \lambda \right) }x.$ Then, for any $\epsilon
>0, $ $y$ has a representation $y=\left\{ y_{k}\right\} _{k\in \mathbb{Z}}$%
such that%
\begin{equation*}
\dsum_{k}\frac{J\left( \lambda ^{k},y_{k}\right) }{K\left( \lambda
^{k},x\right) }<1+\epsilon .
\end{equation*}%
By Hahn-Banach's theorem choose $l_{k}\in \left( \Sigma \overline{X}\right)
^{^{\prime }}$such that $l_{k}\left( x\right) =1,\left\vert l_{k}\left(
a\right) \right\vert \leq \frac{K\left( \lambda ^{k},a\right) }{K\left(
\lambda ^{k},x\right) }$ for $a\in \Sigma \overline{X}.$ In particular we
then have $\left\Vert l_{k}\right\Vert _{X_{i}^{^{\prime }}}\leq \frac{%
\lambda ^{ki}}{K\left( \lambda ^{k},x\right) }$ $\left( i=0,1\right) .$
Define the operator $T$ by the formula $Ta=\dsum_{k}l_{k}\left( a\right)
y_{k}.$ Then 
\begin{eqnarray*}
\dsum_{k}\left\Vert l_{k}\right\Vert _{X_{i}^{^{\prime }}}\left\Vert
y_{k}\right\Vert _{Y_{i}} &\leq &\dsum_{k}\frac{\max \lambda ^{ki}\left\Vert
y_{k}\right\Vert _{Y_{i}}}{K\left( \lambda ^{k},x\right) } \\
&=&\dsum_{k}\frac{J\left( \lambda ^{k},y_{k}\right) }{K\left( \lambda
^{k},x\right) }<1+\epsilon
\end{eqnarray*}%
so $T:\overline{X}\overset{n}{\rightarrow }\overline{y}$ with $\left\Vert
T\right\Vert _{n}<1+\epsilon .$ Since clearly $Tx=y$ (in view of $%
l_{k}\left( x\right) =1$ and $y=\dsum_{k}y_{k}),$moreover $\epsilon >0$
being arbitrary, we conclude that indeed $y<_{n}x.$ The proof is complete. \#

$\medskip $

In particular we get the following corollaries.

\underline{Corollary 1}. $y<<_{n}x$ iff $y<<_{J/K\left( \lambda \right) }x.$

\underline{Corollary 2}. The analogue of \underline{$\frac{1}{2}.2$}$\left(
1\right) $ with $<_{n}$ in place of $<_{J/K\left( \lambda \right) }.\medskip 
$

\underline{Remark}. Bergh in his unfortunately unpublished master's thesis 
\cite{2-BE71} (cf. \cite{3-BeL0}) showed that if for some $t$ holds $J\left(
t,y\right) \leq K\left( t,x\right) $ then $y<_{b}x.$He gave also a number of
interesting concrete applications of this result. Bergh's result apparently
is a special case of the above lemma. With $t=\lambda ^{k}$ apply it to the
trivial representation $y$ such that $y_{j}=\left\{ 
\begin{array}{c}
y:j=k \\ 
0:j\neq k%
\end{array}%
\right\} .$ Cwikel in \cite{5-Cw} on the other hand proved that if we have $%
K\left( \lambda ^{k},y\right) \leq c_{k}K\left( \lambda ^{k},x\right) $ with 
$\dsum_{k}c_{k}<\infty $ then $y<_{b}x.\footnote{%
See Notes \underline{4}:$\left\langle 1\right\rangle $}$ This follow from
our result by taking $\widehat{y}$ to be the representation provided by the
fundamental lemma (\cite{3-BeL0}, p.33). The results of Bergh and Cwikel now
gain a new dimension in the light of our notion of nuclearity.$\medskip $

\underline{4.2}. We can now prove our principal result.$\medskip $

\underline{Theorem}. The couples $\overline{X}$ and $\overline{Y}$ satisfies
condition $\left( O\right) $ iff $y<<_{b}x$ implies $y<<_{n}x.$

\underline{Proof}: Let the couples $\overline{X}$ and $\overline{Y}$ satisfy
condition $\left( O\right) .$Suppose that $y<<_{b}x.$Then by \underline{3.2}%
, prop. we have $y<<_{J/K\left( \lambda \right) }x.$Therefore $y<<_{n}x$ by 
\underline{4.1}, cor. 1. Thus $y<<_{b}x$ implies $y<<_{n}x.$

Conversly assume that this is the case. Then reversing the previous
reasoning we find that $\overline{X}$ and $\overline{Y}$ indeed satisfy
condition $\left( O\right) .$\#$\medskip $

\underline{4.3}. We apply our previous result to give a proof of
Ovchinnikov's theorem \cite{23-Ovc}; until now we have thus assumed this
theorem to be known. Our proof is just a variations of Ovchinnikov's proof
but perhaps slightly simpler; in contract to Janson \cite{16-Jan} it still
involves Grothendieck's fundamental theorem \cite{14-Gr}. $\medskip $

In sec. \underline{3} we have seen that it suffices yo consider the case of $%
\lambda $ -adic weights, i.e. the case $\overline{w}=\left( 1,\lambda
^{-k}\right) .$ We thus have to show that the couples $\overline{l}_{\lambda
}^{\infty }$ and $\overline{l}_{\lambda }^{1}$ satisfy condition $\left(
O\right) .$ Let $y<<_{b}x$. i.e. for some operator $T:\overline{l}_{\lambda
}^{\infty }\rightarrow \overline{l}_{\lambda }^{1}$ holds $y=Tx$ where $x$
and $y$ are given elements in $\Sigma \overline{l}_{\lambda }^{\infty }$ and 
$\Sigma \overline{l}_{\lambda }^{1}$ respectively. We are going to show that 
$y<_{J/K\left( \lambda \right) }x.$ By \underline{3}, prop. it will then
follows that $\overline{l}_{\lambda }^{\infty }$ and $\overline{l}_{\lambda
}^{1}$ indeed satisfy condition $\left( O\right) .$

To achieve this let $\overline{H}=\left( H_{0},H_{1}\right) $ be a suitable
Hilbert couple (i.e. $H_{0}$ and $H_{1}$ are Hilbert spaces). Pick up an
element $u\in \Sigma \overline{H}$ such that $x<<_{b}u$ and $u<<_{K}x.$ That
such a $u$ exists follows at one from the "universal property" of the couple 
$\overline{l}_{\lambda }^{\infty }$ (see \cite{27-Pe}, \cite{28-Pe}, \cite%
{6-CwPe}; cf. Ovchinnikov's treatment \cite{23-Ovc}). In particular we have
then also a bounded linear operator $S:\overline{H}\rightarrow \overline{l}%
_{\lambda }^{\infty }$ with $x=Su.$ The composition thus has the
factorization $TS:\overline{H}\rightarrow \overline{l}_{\lambda }^{\infty
}\rightarrow \overline{l}_{\lambda }^{1}.$ Especially we have $%
TS:H_{i}\rightarrow l^{\infty }\left( \lambda ^{-ik}\right) \rightarrow
l^{1}\left( \lambda ^{-ik}\right) $ $\left( i=0,1\right) .$ Grothendieck's
fundamental theorem has as consequence (\cite{14-Gr}, p. 65\footnote{%
See Notes \underline{5}:$\left\langle 2\right\rangle $}) that $%
TS:H_{i}\rightarrow l^{1}\left( \lambda ^{-ik}\right) $ is indeed nuclear.
Thus $TS:\overline{H}\rightarrow \overline{l}_{\lambda }^{1}$ is "separately
nuclear". In sec. \underline{$\frac{1}{2}$} we said that in general
"separately nuclear" does not entail "nuclear". But in this special case it
does.\footnote{%
See Notes \underline{4}:$\left\langle 3\right\rangle $}. So we have
established $TS:\overline{H}\overset{n}{\rightarrow }\overline{l}_{\lambda
}^{1}.$ It follows that $y<<_{n}u.$ Since ~$u<<_{K\left( \lambda \right) }x$
the ideal property for the ordering $<<_{n}$ (sec. \underline{4.1} cor. 2)
gives $y<<_{J/K\left( \lambda \right) }x.$ The proof is complete.

\section*{5. Janson's proof of Ovchinnikov's theorem.}

We have already given one proof of Ovchinnikov's theorem \cite{23-Ovc},
essentially his own proof; see \underline{4.3}. For the benefit of mankind
we give here another one, essentially the one of Janson \cite{16-Jan}. It is
entirely self-contained in particular independent of the Aronszajn-Gagliardo
theorem \cite{1-AG65}. In fact we can without any extra labor directly
threat the continuous case of sec. \underline{2.4}. Thus $L^{\infty }\left( 
\overline{w}D\right) $ and $L^{1}\left( \overline{z}E\right) $ having the
same meaning as there our goal is to establish the following result.$%
\medskip $

\underline{Theorem}. If $T:L^{\infty }\left( \overline{w}D\right)
\rightarrow L^{1}\left( \overline{z}E\right) $ then $T:L^{\infty }\left(
\rho \left( \overline{w}\right) D\right) \rightarrow L^{1}\left( \rho \left( 
\overline{z}\right) E\right) $ for any $\rho \in \mathcal{P}$ .

$\medskip $

The proof will be broken up into a series of lemmata.$\medskip $

The key lemma, due to Janson, replaces Grothendieck's fundamental theorem in
Ovchinnikov's treatment.

\underline{Lemma 1}. (\cite{16-Jan}, lemma 8\footnote{%
See Notes \underline{5}:$\left\langle 1\right\rangle $}). Let $T:L^{\infty
}\left( \overline{w}\right) \rightarrow l^{1}\left( \overline{z}\right) .$%
Then $T:l^{\infty }\rightarrow l^{1}$ provided the following condition is
fulfilled:

$\left( 1\right) :\dsum_{m\in \mathbb{Z}}\epsilon _{m}<\infty $ for some
(positive) sequence $\left\{ \epsilon _{m}\right\} $ such that $\min \left(
w_{k}^{0}/z_{j}^{0},w_{k}^{1}/z_{j}^{1}\right) \leq \epsilon _{m}$ if $%
j-k=m.\medskip $

\underline{Remark}. As will be clear from the proof the same conclusion
holds if we just assume $T:c_{0}\left( \overline{w}\right) \rightarrow
l^{1}\left( \overline{z}\right) $ ($c_{0}$ denotes of course the space of
(doubly infinite) sequences tending to $0.).$We see also that $T:l^{\infty
}\rightarrow l^{1}$ is in fact a \underline{nuclear} operator; this
observation is perhaps new.$\medskip $

\underline{Proof} (after Janson \cite{16-Jan}): On the space $\Sigma
c_{0}\left( \overline{w}\right) =c_{0}\left( \min \left( w_{0},w_{1}\right)
\right) $ the operator $T$ is given by a (doubly infinite) matrix $\left(
t_{jk}\right) ,$so that if $\left\{ x_{j}\right\} \in \Sigma \left(
c_{0}\left( \overline{w}\right) \right) $ we have $Tx_{j}=%
\dsum_{k}t_{jk}x_{k}$ which again entails that 
\begin{equation}
\dsum_{j}z_{j}^{i}\left\vert \dsum_{k}t_{jk}x_{k}\right\vert \leq
C\sup_{k}w_{k}^{i}\left\vert x_{k}\right\vert ,\left( i=0,1\right)   \tag{2}
\end{equation}%
for some constant $C.$Taking $j=0$ in $\left( 1\right) $ we see in particular%
\begin{equation}
\min \left( w_{k}^{0},w_{k}^{1}\right) \rightarrow 0\text{ as }\left\vert
k\right\vert \rightarrow \infty .  \tag{$1^{^{\prime }}$}
\end{equation}%
Therefore we have $l^{\infty }\subseteq \Sigma \left( c_{0}\left( \overline{w%
}\right) \right) .$To complete the proof it thus suffices to show that $%
\dsum_{j,k}\left\vert t_{jk}\right\vert <\infty ,$ which again - in view of $%
\left( 1\right) $ - will follows if we can can show that%
\begin{equation}
\dsum_{j-k=m}\left\vert t_{jk}\right\vert \leq C\epsilon _{m}  \tag{3}
\end{equation}%
To this end we set $a_{jk}=\max \left(
z_{j}^{0}/w_{k}^{0},z_{j}^{1}/w_{k}^{1}\right) t_{jk}.$Then $\left( 2\right) 
$ can be written simply as 
\begin{equation}
\dsum_{j}\left\vert \dsum_{k}a_{jk}\zeta _{k}\right\vert \leq C\sup
\left\vert \zeta _{k}\right\vert .  \tag{$2^{^{\prime }}$}
\end{equation}%
That is, the matrix $\left( a_{jk}\right) $ defines a bounded linear
operator $A:l^{\infty }\rightarrow l^{1}.$We see that in order to establish $%
\left( 3\right) $ it suffices to show that $\left( 2^{^{\prime }}\right) $
implies%
\begin{equation}
\dsum_{j-k=m}\left\vert a_{jk}\right\vert \leq C  \tag{$3^{^{\prime }}$}
\end{equation}%
for each $m.$ Now $\left( 2^{^{\prime }}\right) $ entails that for each real
number $x$ holds%
\begin{equation*}
\dsum_{j}\left\vert \dsum_{k}a_{jk}e^{i\left( j-k-m\right) x}\right\vert
\leq C
\end{equation*}%
Taking averages $\left( 3^{^{\prime }}\right) $ then follows at once. \#$%
\medskip $

\underline{Remark}. As already mentioned (see Introduction) a different
"elementary" method in connection with Ovchinnikov's theorem \cite{23-Ovc}
based on Khinchine's (=Xincin's) inequality is used in Gustavsson \cite%
{15-Gu}.$\medskip $

Now we restate lemma 1 in a form more useful for the application we have in
mind.

\underline{Lemma $1^{^{\prime }}$}. Let $T:\left( l^{\infty },l^{\infty
}\left( \frac{1}{\tau }\right) \right) \rightarrow \left( l^{1},l^{1}\left( 
\frac{1}{\sigma }\right) \right) .$Then $T:l^{\infty }\left( \frac{1}{\rho
\left( \tau \right) }\right) \rightarrow l^{1}\left( \frac{1}{\lambda \left(
\sigma \right) }\right) $ provided $\dsum_{m\in \mathbb{Z}}\epsilon
_{m}<\infty $ where 
\begin{equation*}
\min \left( 1,\frac{\sigma _{j}}{\tau _{k}}\right) \frac{\rho \left( \tau
_{k}\right) }{\lambda \left( \sigma _{j}\right) }\leq \epsilon _{m}\text{ if 
}k-j=m
\end{equation*}%
Here $\tau =\left\{ \tau _{j}\right\} $ and $\sigma =\left\{ \sigma
_{j}\right\} $ are any two given sequences and $\rho $ and $\lambda $ are
given positive functions on $\left( 0,\infty \right) ;\rho \left( \tau
\right) $ stands for the sequence $\left\{ \rho \left( \tau _{k}\right)
\right\} $ with a similar meaning for $\lambda \left( \sigma \right) ;1$
denotes the sequence $\left\{ 1\right\} ,\frac{1}{\tau }$ the sequence $%
\left\{ \frac{1}{\tau _{k}}\right\} ,$ similarly for $\frac{1}{\sigma },%
\frac{1}{\rho \left( \lambda \right) },\frac{1}{\lambda \left( \sigma
\right) }.$

\underline{Proof}: Indeed we at once reduce it to lemma $1$ with $w^{0}=\rho
\left( \tau \right) ,w^{1}=\frac{\rho \left( \tau \right) }{\tau }%
,z^{0}=\lambda \left( \sigma \right) ,z^{1}=\frac{\lambda \left( \sigma
\right) }{\sigma }.\#\medskip $

In particular if $\tau =\sigma ,\rho =\lambda $ condition $\left( 1\right) $
takes the form%
\begin{equation}
\dsum_{m}\epsilon _{m}<\infty \text{ where }\min \left( 1,\frac{\tau _{j}}{%
\tau _{k}}\right) \frac{\rho \left( \tau _{k}\right) }{\rho \left( \tau
_{j}\right) }\text{ if }k-j=m.  \tag{$1^{^{\prime \prime }}$}
\end{equation}%
$\medskip $

\underline{Example}. Notice that $\left( 1^{^{\prime \prime }}\right) $
certainly is fulfilled if $\tau _{k}=\lambda ^{k},\rho \in \mathcal{P}^{+-}.$
Indeed denoting (for any $\rho )$ by $s_{\rho }$ the corresponding dilation
function \footnote{%
See Notes \underline{5}.$\left\langle 2\right\rangle $} i.e.%
\begin{equation}
s_{\rho }\left( t\right) =\sup_{u}\frac{\rho \left( tu\right) }{\rho \left(
u\right) },
\end{equation}%
we have $\min \left( 1,\lambda ^{j-k}\right) \frac{\rho \left( \lambda
^{k}\right) }{\rho \left( \lambda ^{j}\right) }\leq \min \left( 1,\lambda
^{-m}\right) s_{\rho }\left( \lambda ^{m}\right) $ for $k-j=m$ so $\left(
1^{^{\prime \prime }}\right) $ follows if $\rho \in \mathcal{P}^{+-}$
because then $s_{\rho }\left( t\right) =O\left( \min \left( t^{\alpha
_{0}},t^{\alpha _{1}}\right) \right) $ for suitable $\alpha _{0},\alpha _{1}$
with $0<\alpha _{0}<\alpha _{1}<1.\medskip $

\underline{Lemma 2}. Let $x\in \Sigma L^{\infty }\left( \overline{w}D\right) 
$ and let $\tau =\left\{ \tau _{k}\right\} $ be any \underline{increasing}
positive sequence such that%
\begin{equation}
K\left( t,x\right) \leq C\sup_{k}\min \left( 1,\frac{t}{\tau _{k}}\right)
K\left( \tau _{k},x\right)  \tag{4}
\end{equation}%
for some $C.$Then we can find an element $a\in \Sigma l^{\infty }\left( 1,%
\frac{1}{\tau }\right) $ and a bounded linear map $\iota :l^{\infty }\left(
1,\frac{1}{\tau }\right) \rightarrow L^{\infty }\left( \overline{w}D\right) $
with $x=\iota a$ such that $a<<_{K}x.$

\underline{Proof}: We first assume that the underlying measure space is $%
W=\left( 0,\infty \right) $ and that the measure is $\mu =dt$ and take $D$
to be the trivial bundle all of whose fibers are equal to $\mathbb{R}$. We
thus have the couple $L^{\infty }\left( 1,\frac{1}{\tau }\right) .$We
further take $w^{0}=1,w^{1}=\frac{1}{t}.$ Then we have the well-known
estimates (see e.g. \cite{30-Pe})%
\begin{equation}
K\left( t,x\right) \approx \sup_{s\in \left( 0,\infty \right) }\min \left( 1,%
\frac{t}{s}\right) \left\vert x\left( s\right) \right\vert \text{ for }x\in
\Sigma L^{\infty }\left( 1,\frac{1}{t}\right) .  \tag{5}
\end{equation}%
We further assume that $x\in \mathcal{P}_{1},$ which in view of $\left(
5\right) $ implies that $K\left( t,x\right) \approx x\left( t\right) $ and
also \underline{linear} in the intervals $\left( \tau _{k},\tau
_{k+1}\right) .$(Then $\left( 4\right) $ is automatically fulfilled). We set 
$a=\left\{ x\left( \tau _{k}\right) \right\} $ and define the map $\iota $
as follows. If $b=\left\{ b_{k}\right\} \in \Sigma l^{\infty }\left( 1,\frac{%
1}{\tau }\right) $ we put $\iota b\left( \tau _{k}\right) =b_{k}$ and extend
this by linearity to general $t,$ that is if $t$ is in the interval $\left(
\tau _{k},\tau _{k+1}\right) $ we put $\iota b\left( t\right) =\left(
1-\theta \right) b_{k}+\theta b_{k+1}$ where we have written $t=\left(
1-\theta \right) \tau _{k}+\theta \tau _{k+1}$ with $0<\theta <1.$It is
clear that%
\begin{eqnarray*}
\left\vert \iota b\left( t\right) \right\vert  &\leq &\left( 1-\theta
\right) \left\vert b_{k}\right\vert +\theta \left\vert b_{k}\right\vert \leq
\sup_{k}\left\vert b_{k}\right\vert =\left\Vert b\right\Vert _{l^{\infty },}
\\
\frac{1}{t}\left\vert \iota b\left( t\right) \right\vert  &\leq &\left(
1-\theta \right) \frac{\tau _{k}}{t}\frac{\left\vert b_{k}\right\vert }{\tau
_{k}}+\theta \frac{\tau _{k+1}}{t}\frac{\left\vert b_{k+1}\right\vert }{%
t_{k+1}}\leq  \\
&\leq &\sup_{k}\frac{\left\vert b_{k}\right\vert }{\tau _{k}}=\left\Vert
b\right\Vert _{l^{\infty }\left( \frac{1}{\tau }\right) }
\end{eqnarray*}%
which establishes the continuity of $\iota $, and that $\iota a=x.$ Also by
the discrete analogue of $\left( 5\right) $ and by $\left( 4\right) $%
\begin{equation*}
K\left( t,a\right) \approx \sup_{k}\min \left( 1,\frac{t}{\tau _{k}}\right)
x\left( \tau _{k}\right) \approx x\left( t\right) \approx K\left( t,x\right) 
\end{equation*}%
so in particular $a<<_{K}x.$

To treat the general case we shall make use of the fact that the couples $%
L^{\infty }\left( 1,\frac{1}{t}\right) =\left( L^{\infty },L^{\infty }\left( 
\frac{1}{t}\right) \right) $ and $L^{\infty }\left( \overline{w}D\right) $
are relative Calderon. In the case of a trivial bundle with all fibers equal
to $\mathbb{R}$ ("scalar case", this is well-known (see \cite{27-Pe}, \cite%
{28-Pe}, \cite{6-CwPe}). The general case can easily brought back to this
case, since $L^{\infty }\left( \overline{w}D\right) $ is a partial retract
of a scalar couple (cf. \underline{2.4}).

Let thus $x$ be given in $\Sigma L^{\infty }\left( \overline{w}D\right) .$%
Let $\widetilde{x}$ be the function in $\Sigma L^{\infty }\left( 1,\frac{1}{t%
}\right) $ which is linear in each interval $\left( \tau _{k},\tau
_{k+1}\right) $ and such that $\widetilde{x}\left( \tau _{k}\right) =K\left(
\tau _{k},x\right) .$ It is clear from $\left( 4\right) $ that $x<<_{K}%
\widetilde{x}$ so we get a map $T:L^{\infty }\left( 1,\frac{1}{t}\right)
\rightarrow L^{\infty }\left( \overline{w}D\right) $ with $x=T\widetilde{x}.$
Also trivially $\widetilde{x}<<_{K}x.$Furthermore the previous proof on the
special case $Z=\left( 0,\infty \right) $ yields a map $\widetilde{\iota }%
:l^{\infty }\left( 1,\frac{1}{\tau }\right) \rightarrow L^{\infty }\left( 1,%
\frac{1}{t}\right) $ with $\widetilde{x}=\widetilde{\iota }a$ for some $a\in
\Sigma \left( l^{\infty },l^{\infty }\left( \frac{1}{\tau }\right) \right) $
with $a<<_{K}\widetilde{x}.$ We obtain the desired map $\iota :l^{\infty
}\left( 1,\frac{1}{\tau }\right) \rightarrow L^{\infty }\left( \overline{w}%
D\right) $ as the composite $T\widetilde{\iota }.$The proof is completed if
we notice that by the transitivity of the ordering $<<_{K}$ (sec. $%
\underline{\frac{1}{2}.2}$ ) follows $a<<_{K}x.$\#$\medskip $

The result dual to lemma 2 reads.

\underline{Lemma 3}.Let $y\in \Sigma L^{1}\left( \overline{z}E\right) $ and
let $\tau =\left\{ \tau _{k}\right\} $ be any \underline{increasing}
positive sequence such that%
\begin{equation}
\inf_{k}\max \left( 1,\frac{t}{\tau _{k}}\right) \rho \left( \tau
_{k}\right) \leq C\rho \left( t\right)  \tag{6}
\end{equation}%
for some constant $C.$Then we can find an element $b\in \Sigma l^{1}\left( 1,%
\frac{1}{\tau }\right) $ and a bounded linear map $\pi :L^{1}\left( 
\overline{z}E\right) \rightarrow l^{1}\left( 1,\frac{1}{\tau }\right) $ with 
$b=\pi y$ such that if $b\in l^{1}\left( \frac{1}{\rho \left( \tau \right) }%
\right) $ then $y\in L^{1}\left( \rho \left( \overline{z}\right) E\right) .$

\underline{Proof}: Let $\left\{ H_{k}\right\} $ be a partition of the
underlying measure space $Z$ into disjoint subsets such that%
\begin{equation}
\rho \left( \tau _{k}\right) \max \left( 1,\frac{z^{0}}{z^{1}\tau _{k}}%
\right) \leq C\rho \left( \frac{z^{0}}{z^{1}}\right)  \tag{7}
\end{equation}%
on $H_{k}.$ (It is for the existence of such a partition we need $\left(
6\right) ).$

By the Hahn-Banach theorem (cf. 2.2, proof of prop. 1) we get linear
functionals. $\beta _{k}$ $\left( k\in \mathbb{Z}\right) $ on $\Sigma
L^{1}\left( \overline{z}E\right) $ such that%
\begin{equation*}
\beta _{k}\left( y\right) =\tint\nolimits_{H_{k}}z^{0}/\rho \left(
z^{0}/z^{1}\right) \left\Vert y\right\Vert _{E_{\omega }}dv
\end{equation*}%
and%
\begin{equation*}
\left\vert \beta _{k}\left( v\right) \right\vert \leq
\tint\nolimits_{H_{k}}z^{0}/\rho \left( z^{0}/z^{1}\right) \left\Vert
v\right\Vert _{E_{\omega }}dv
\end{equation*}%
for every $v\in \Sigma L^{1}\left( \overline{z}E\right) .$ Then by $\left(
7\right) $ we have also 
\begin{eqnarray*}
\left\vert \beta _{k}\left( v\right) \right\vert &\leq &\frac{C}{\rho \left(
\tau _{k}\right) }\tint\nolimits_{H_{k}}z^{0}\left\Vert v\right\Vert
_{E_{\omega }}dv \\
\left\vert \beta _{k}\left( v\right) \right\vert &\leq &\frac{C\tau _{k}}{%
\rho \left( \tau _{k}\right) }\tint\nolimits_{H_{k}}z^{1}\left\Vert
v\right\Vert _{E_{\omega }}dv
\end{eqnarray*}%
We define $\pi $ by the formula $\pi v_{k}=\rho \left( \tau _{k}\right)
\beta _{k}\left( v\right) $ for $v\in \Sigma L^{1}\left( \overline{z}%
E\right) .$Then obviously 
\begin{eqnarray*}
\left\Vert \pi v\right\Vert _{l^{1}} &=&\dsum_{k}\left\vert \pi
v_{k}\right\vert \leq C\dsum_{k}\tint\nolimits_{H_{k}}z^{0}\left\Vert
v\right\Vert _{E_{\omega }}dv=C\left\Vert v\right\Vert _{L^{1}\left(
z^{0}E\right) } \\
\left\Vert \pi v\right\Vert _{l^{1}\left( \frac{1}{\tau }\right) }
&=&\dsum_{k}\frac{1}{\tau _{k}}\left\vert \pi v_{k}\right\vert \leq
C\dsum_{k}\tint\nolimits_{H_{k}}z^{1}\left\Vert v\right\Vert _{E_{\omega
}}dv=C\left\Vert v\right\Vert _{L^{1}\left( z^{1}E\right) }
\end{eqnarray*}%
Moreover for $\pi y=b=\left\{ b_{k}\right\} $ we obtain%
\begin{eqnarray*}
\left\Vert b\right\Vert _{l^{1}\left( \frac{1}{\rho \left( \tau \right) }%
\right) } &=&\dsum_{k}\frac{\left\vert b_{k}\right\vert }{\rho \left( \tau
_{k}\right) }=\dsum_{k}\tint\nolimits_{H_{k}}z^{0}/\rho \left(
z^{0}/z^{1}\right) \left\Vert y\right\Vert _{E_{\omega }}dv \\
&=&\tint\nolimits_{Z}z^{0}/\rho \left( z^{0}/z^{1}\right) \left\Vert
y\right\Vert _{E_{\omega }}dv=\left\Vert y\right\Vert _{L^{1}\left( \rho
\left( \overline{z}\right) E\right) }
\end{eqnarray*}%
Thus $y\in L^{1}\left( z^{0}/\rho \left( z^{0}/z^{1}\right) E\right) $ iff $%
b\in l^{1}\left( \frac{1}{\rho \left( \tau \right) }\right) .$\#$\medskip $

Finally we establish the following result also taken over in principle from
Janson \cite{16-Jan}.

\underline{Lemma 4}. Let $x$ be in $\Sigma L^{\infty }\left( \overline{w}%
D\right) $ and set $\rho \left( t\right) =K\left( t,x\right) $ (so $\rho \in 
\mathcal{P)}$. Assume the range of the function $\rho \left( t\right) $ as
well as that of $\rho ^{^{\prime }}\left( t\right) =\frac{t}{\rho \left(
t\right) }$ is $\left( 0,\infty \right) .$Then there exists an \underline{%
increasing} positive sequence $\tau =\left\{ \tau _{k}\right\} _{k\in 
\mathbb{Z}}$ such that $\left( 4\right) ,\left( 6\right) $ and $\left(
1^{^{\prime \prime }}\right) $ are fulfilled.

\underline{Proof}: We begin by noticing that $\left( 4\right) $ says that%
\begin{equation}
\rho \left( t\right) \leq C\cdot \inf_{k}\min \left( 1,\frac{t}{\tau _{k}}%
\right) \rho \left( \tau _{k}\right)  \tag{8}
\end{equation}%
and that $\left( 6\right) $ is formally $\left( 8\right) $ for the function $%
\rho ^{^{\prime }}.$ Following Janson \cite{16-Jan} we define $\tau _{k}$
inductively by $\tau _{0}=1$ and $\min \left( \frac{\rho \left( \tau
_{k+1}\right) }{\rho \left( \tau _{k}\right) },\frac{\rho ^{^{\prime
}}\left( \tau _{k+1}\right) }{\rho ^{^{\prime }}\left( \tau _{k}\right) }%
\right) =2.$ (It is here the assumption about the range comes in.). Then%
\begin{equation*}
\rho \left( \tau _{j}\right) \geq 2^{j-k}\rho \left( \tau _{k}\right) \text{
or }\frac{\rho \left( \tau _{k}\right) }{\rho \left( \tau _{j}\right) }\leq
2^{k-j}\text{ if }j\geq k
\end{equation*}%
Similarly%
\begin{equation*}
\rho ^{^{\prime }}\left( \tau _{k}\right) \geq 2^{k-j}\rho ^{^{\prime
}}\left( \tau _{j}\right) \text{ or }\frac{\tau _{j}\rho \left( \tau
_{k}\right) }{\tau _{k}\rho \left( \tau _{j}\right) }\leq 2^{j-k}\text{ if }%
j\leq k
\end{equation*}%
Thus at any rate%
\begin{equation*}
\min \left( 1,\frac{\tau _{j}}{\tau _{k}}\right) \frac{\rho \left( \tau
_{k}\right) }{\rho \left( \tau _{j}\right) }\leq 2^{-\left\vert
j-k\right\vert }
\end{equation*}%
so $\left( 1^{^{\prime \prime }}\right) $ follows indeed, with $\epsilon
_{m}=2^{-\left\vert m\right\vert }.$To establish $\left( 6\right) $ let $%
t\in \left( \tau _{k},\tau _{k+1}\right) .$Then there are two cases:

\underline{Case 1}. $\rho \left( \tau _{k+1}\right) \leq 2\rho \left( \tau
_{k}\right) .$Then, since $\rho \in \mathcal{P}_{1},$ we get $\rho \left(
t\right) \leq \rho \left( \tau _{k+1}\right) \leq 2\rho \left( \tau
_{k}\right) =2\min \left( 1,\frac{t}{\tau _{k}}\right) \rho \left( \tau
_{k}\right) .$

\underline{Case 2}. $\frac{\rho \left( t_{k+1}\right) }{\tau _{k+1}}\leq
2\cdot \frac{\rho \left( \tau _{k}\right) }{\tau _{k}}.$We now have $\frac{%
\rho \left( t\right) }{t}\leq 2\frac{\rho \left( \tau _{k+1}\right) }{\tau
_{k+1}}$ or $\rho \left( t\right) \leq 2\frac{t}{\tau _{k+1}}\rho \left(
\tau _{k+1}\right) =2\min \left( 1,\frac{t}{\tau _{k+1}}\right) \rho \left(
\tau _{k+1}\right) .$

Therefore we have proven $\left( 8\right) $ (or $\left( 4\right) )$ with $%
C=2.$Since $\left( 6\right) ,$ as we have just remarked, is just $\left(
8\right) $ with $\rho ^{^{\prime }}$ in place of $\rho $ and further $\rho $
and $\rho ^{^{\prime }}$ enter in a symmetric fashion into the definition of 
$\tau $ we need no particular proof of $\left( 6\right) .$ \#$\medskip $

These lemmas being established we can at last proceed to the

\underline{Proof} (of theorem/completed/). Let thus $T:L^{\infty }\left( 
\overline{w}D\right) \rightarrow L^{1}\left( \overline{z}E\right) $ and $%
x\in L^{\infty }\left( \rho \left( \overline{w}\right) D\right) .$ We have
to show that $Tx\in L^{1}\left( \rho \left( \overline{z}\right) E\right) .$
Clearly it suffices to prove this if $\rho \left( t\right) =K\left(
t,x\right) .$(cf. e.g. the argument in \underline{1}, proof of prop. 2).
Make first the assumption about the ranges in lemma 4. Let then $\tau $ be
the sequence provided by that lemma and select according to lemma 2, $b,\pi $
according to lemma $3,$ with $y=Tx.$ Put $S=\pi T$ so that we have the
commutative diagram%
\begin{equation*}
\begin{tabular}{ccc}
$L^{\infty }\left( \overline{w}D\right) $ & $\overset{T}{\rightarrow }$ & $%
L^{1}\left( \overline{z}E\right) $ \\ 
$\uparrow \iota $ &  & $\downarrow \pi $ \\ 
$l^{\infty }\left( 1,\frac{1}{\tau }\right) $ & $\overset{S}{\rightarrow }$
& $l^{1}\left( 1,\frac{1}{\tau }\right) $%
\end{tabular}%
\end{equation*}%
By lemma $1$ $S:l^{\infty }\left( \frac{1}{\rho \left( \tau \right) }\right)
\rightarrow l^{1}\left( \frac{1}{\rho \left( \tau \right) }\right) .$Now $%
a\in l^{\infty }\left( \frac{1}{\rho \left( \tau \right) }\right) $ so $%
Sa\in l^{1}\left( \frac{1}{\rho \left( \tau \right) }\right) .$But $Sa=\pi
T\iota a=\pi Tx=\pi y.$ Therefore $y\in L^{1}\left( \rho \left( \overline{z}%
\right) E\right) .$

This completes the proof in the above assumption about the ranges. If this
assumptions is \underline{not} fulfilled one must take recourse to suitable
one-sided version of the preceding lemmata. In order to avoid tedious
repetition we omit the details.$\#$

\section*{6. Open questions.}

\underline{6.1}. \underline{"The philosophy of weights"}. We would like to
give an at least heuristic explanation why weights play such a great role in
the theory of interpolation spaces. We proceed in several steps.$\medskip $

Q: Which are the simplest pairs?

A: The pairs $\overline{A}=\left( A_{0},A_{1}\right) $ for which $A_{0}$ and 
$A_{1}$ are similar. By the latter we mean that there exists a Banach space $%
E$ and linear maps $\Lambda _{i}:A_{i}\rightarrow E$ $\left( i=0,1\right) $
such that $\left\Vert a\right\Vert _{A_{i}}=\left\Vert \Lambda
_{i}a\right\Vert _{E},a\in A_{i},$ i.e. $A_{i}$ is isometric to a (closed)
subspace of $E.\medskip $

\underline{Example}. If the $A_{i}$ are Hilbert spaces then this is always
the case (with $E\ $Hilbert too); cf. \cite{11-FoLi}.

In this context it is interesting to consider interpolation spaces $A$ with
respect to $\overline{A}$ such that $\left\Vert a\right\Vert _{A}=\left\Vert
\Lambda a\right\Vert _{E},a\in A,$ for a suitable linear map $\Lambda .$
Indeed $\Lambda $ may be conceived as a kind of function of $\Lambda _{0}$
and $\Lambda _{1}.$ We have then the \underline{problem of interpolation
functions}, a line of development starting with the classical paper by
Foias-Lions \cite{11-FoLi} already referred to.$\medskip $

Q. Which are the simplest Banach spaces? And which are the simplest
operators?

A. $E=l^{1},l^{\infty },l^{2},l^{p}$ and there continuous analogues $%
L^{1},L^{\infty },L^{2},L^{p}$, more generally rearrangement invariant
spaces (both sequence and function spaces), in particular thus Marcinkiewicz
and Lorentz spaces, $BMO,H^{1},FL^{1},FL^{\infty }$ etc. As for $\Lambda $ a
natural candidate in each of these cases is the multiplication operator $%
\Lambda f=wf.\medskip $

We have thus a general program, albeit somewhat vague in its contours. 
\underline{To generalize the } \underline{theory of Ovchinnikov}  \underline{%
and Janson to the context of more general weighted pairs} $E\left( \overline{%
w}\right) .\medskip $

Let us consider a few somewhat more specialized questions.

Is $E\left( \overline{w}\right) $ always Calderon?. And when is $E\left( 
\overline{w}\right) $ tame in Ovchinnikov's sense \cite{23-Ovc}. (The case $%
FL^{1}$ treated by Janson \cite{16-Jan} shows that $E\left( \overline{w}%
\right) $ is not always tame.)

To relate the properties of $E\left( \overline{w}\right) $ with the dual
pair $E^{\ast }\left( 1/\overline{w}\right) .$

To characterize all (relative) interpolation spaces with respect to the
pairs $l^{p}\left( \overline{w}\right) $ and $l^{q}\left( \overline{z}%
\right) .$ (Some results in this direction are given in Ovchinnikov \cite%
{24-Ovc}. His paper \cite{25-Ovc} on the other hand might be considered as a
beginning of a study of the (Hilbert) pair $l^{2}\left( \overline{w}\right) $
in this spirit. As we have already noticed the case $E=FL^{1}$ as well as $%
E=FL^{\infty }$ appears in Janson \cite{16-Jan}.)$\medskip $

\underline{Remark}. The special role of $l^{\infty }\left( \overline{w}%
\right) $ and $l^{1}\left( \overline{z}\right) ,$ further accentuated by the
recent work of Brudnyi-Krugljak. \cite{4-BrKr}, is apparent already in \cite%
{26-Pe} (cf. \cite{6-CwPe} too); from the category point of view advocated
there they simply appear as \underline{projective} and \underline{injective}
objects respectively.$\medskip $

\underline{6.2}. \underline{Lion's problem}. Lions \cite{20-Lio} long ago
raised the question whether the space $\overline{X}_{\theta p}$ always
depends \underline{effectively} on their parameters $\theta ,p.$ This is
vaguely related to the topics of sec. \underline{4}. (Cf. \cite{27-Pe} for
quite a different idea.). For the complex spaces $\overline{X}_{\left[
\theta \right] }$ the analogues question was solved by Stafney \cite{32-Sta}%
. His method is however quite general and in the case of the space $%
\overline{X}_{\theta p}$ leads to the following partial result (cf. \cite%
{3-BeL0} p.82, exc. 21 and p.104): $\left( \ast \right) $ \underline{Assume
there exists} $\theta _{0},p_{0}$ \underline{and} $\theta _{1},p_{1}$ 
\underline{with} $\theta _{0},\theta _{1}$ \underline{in} $\left( 0,1\right)
,\theta _{0}\neq \theta _{1}$ \underline{such that} $\overline{X}_{\theta
_{0}p_{0}}=\overline{X}_{\theta _{1}p_{1}}.$\underline{Assume also that} $%
\overline{X}$ is \underline{regular} . \underline{Then} $X_{0}=X_{1}.$ There
remains thus the more difficult case $\theta _{0}=\theta _{1}$, that is
(writing $\theta =\theta _{0}=\theta _{1})$ the case when $\overline{X}%
_{\theta p_{0}}=\overline{X}_{\theta p_{1}}$for some $\theta $ in $\left(
0,1\right) $ with $p_{0}\neq p_{1}.$ Applying $\left( \ast \right) $ in
conjunction with a suitable reiteration theorem we conclude that $\overline{X%
}_{\theta p}$ in this case is independent of $p$ if $p<\infty ;$ regarding
the limiting case $p=\infty $ we can only say that the norm of $\overline{X}%
_{\theta \infty }$ restricted to $\overline{X}_{\theta 1}$ coincides up to
equivalence with that of $\overline{X}_{\theta 1}.$ On the other hand by a
theorem of Lions (see \cite{3-BeL0}, p. 103) holds 
\begin{eqnarray*}
\left( \overline{X}_{\theta p},\overline{X}_{\theta q}\right) _{\left[ \tau %
\right] } &=&\overline{X}_{\zeta r}\text{ if }\zeta =\left( 1-\tau \right)
\theta +\tau \theta \\
\frac{1}{r} &=&\frac{1-\tau }{p}+\frac{\tau }{q},\left( \tau \in \left(
0,1\right) \right) .
\end{eqnarray*}%
Thus we arrive at least to the following conclusion. If $\overline{X}%
_{\theta p}$ is independent of $p$ for some $\theta $ in $\left( 0,1\right) $
then the same is true for any $\theta $ in $\left( 0,1\right) .\medskip $

If we now generalizing (replacing $\rho \left( t\right) =t^{\theta }$ by a
general function), we have the problem to decide for any particular $\rho
\in \mathcal{P}$ (or for all $\rho \in \mathcal{P)}$ whether one can have $%
\overline{X}_{\rho \infty ;K}=\overline{X}_{\rho 1;J}.$ Let us remark that
if there exists $x\in \Sigma X$ such that $K\left( t,x\right) \approx \rho
\left( t\right) $ then this certainly cannot be the case. (This is because $%
\Delta $ is dense in $\overline{X}_{\rho 1;J}.).$But then we have at least
the modified question whether the norms of $\overline{X}_{\rho \infty ;K}$
and $\overline{X}_{\rho 1;J}$ coincide (up to equivalence) on $\Delta .$ Put
in more concrete terms if $x\in \Sigma \overline{X}$ and $K\left( t,x\right)
=o\left( \rho \left( t\right) \right) $ does there exists a representation $%
\widehat{x}=\left\{ x_{k}\right\} _{k\in \mathbb{\mathbb{Z}}}$ of $x$ such
that%
\begin{equation*}
\dsum_{k}\frac{J\left( 2^{k},x_{k}\right) }{\rho \left( 2^{k}\right) }%
<\infty .
\end{equation*}

This is about to say that the identity map $\overline{X}$ is in some sense
nuclear. Let us remark that if we pass to the more general situation of $%
\underline{\text{locally convex spaces }}$ (instead of just Banach spaces)
this very well can be the case. This is essentially contained already in a
classic example in \cite{18-LiPe}, chap. $VIII.\medskip $

\bigskip \underline{6.3}. We conclude by a few scattered remarks pertaining
mainly to the topics on sec. \underline{3} .

We can ask if $\overline{Y}$ of type $\left( u\right) $ implies $\overline{Y}%
^{^{\prime }}$ (the dual pair, assuming that $\overline{Y}$ is regular) then
is of type $\left( l\right) .\medskip $

We have also a natural notion of "type $\left( l_{0}\right) ^{"}$which is
roughly speaking obtained if we replace \thinspace $l^{\infty }$ by $c_{0}.$
(Janson \cite{16-Jan} makes a great show of the distinction between $%
l^{\infty }$ and $c_{0},$ a distinction which we have here almost entirely
neglected.). And corresponding to it a "condition $\left( O_{0}\right) ^{".}$%
The question is then if $\overline{X}$ is of type $\left( l_{0}\right) $
implies $\overline{X}^{^{\prime }}$ is type $\left( u\right) .$ And if from
the facts that $\overline{X}$ and $\overline{Y}$ satisfies condition $\left(
O_{0}\right) $ follows that $\overline{Y}^{^{\prime }}$ and $\overline{X}%
^{^{\prime }}$ satisfy condition $\left( O\right) .\medskip $

In this connection we also ask if $\overline{X}$ is Calderon is then $%
\overline{X}^{^{\prime }}$ Calderon too.

If $\overline{X}$ is tame in Ovchinnikov's sense \cite{23-Ovc} does it
follows that $G^{A}\left( \overline{X}\right) =H^{B}\left( \overline{X}%
\right) $ where $A$ and $B$ are in the same relation as $l^{\infty }\left(
\rho \left( \overline{w}\right) \right) $ and $l^{1}\left( \rho \left( 
\overline{w}\right) \right) $ (cf. sec. 1).

\section*{Appendix.}

Consider any Banach couple $\overline{A}=\left\{ A_{0},A_{1}\right\} .$It is
well-known (see e.g. \cite{6-CwPe}) $\overline{A}$ is termed \underline{%
Calderon} if $b<<_{K}a$ implies $b<<_{b}a.$ $\overline{A}$ will be called $%
\sigma -$additive if for any (infinite) sequence $\left\{ a_{j}\right\}
\subseteq \Sigma A$ such that for some (or any) $t$ holds $\dsum_{j}K\left(
t,a_{j}\right) <\infty $ we can find a sequence $\left\{ b_{j}\right\}
\subseteq \Sigma A$ such that for any $t$ holds $K\left( t,\dsum
b_{j}\right) \approx \dsum_{j}K\left( t,b_{j}\right) $ and (for each $j)$ $%
K\left( t,a_{j}\right) \approx K\left( t,b_{j}\right) .$\footnote{%
See Notes Appendix:$\left\langle 1\right\rangle .$}. Assuming that $%
\overline{A}$ is both Calderon and $\sigma $-additive one can then easily
prove the following fact (see \cite{3-BeL0}, p. 12 remark 2.2, where the
additive case is considered). If $A$ is any interpolation space with respect
to $\overline{A}$ then there exists a sequence norm $\Phi $ such that 
\begin{equation}
\left\Vert a\right\Vert _{A}\approx \Phi \left( \left\{ K\left( \lambda
^{j},a\right) \right\} \right) .  \tag{1}
\end{equation}%
(Every interpolation is a classical $K$ space.). Suppose moreover that $%
\overline{A}$ satisfies the \underline{strong form of the fundamental lemma}%
, i.e. every $a\in \Sigma _{0}\overline{A}$ has a representation $a=\left\{
a_{k}\right\} $ such that for any $t$ holds 
\begin{equation*}
\dsum \min \left( 1,\frac{t}{\lambda ^{k}}\right) J\left( \lambda
^{k},a_{k}\right) \leq C\cdot K\left( t,a\right) 
\end{equation*}%
with $C$ depending only on $\overline{A}$ (see \cite{6-CwPe}, p. 31, th.
4.4, "property $\left( 2\right) ^{^{\prime \prime }}),$ or $\Omega \left(
\left\{ J\left( \lambda ^{k},a_{k}\right) \right\} \right) \leq C\cdot
K\left( \lambda ^{k},a\right) $ where we have introduced the discrete
Calderon transform $\Omega ,$i.e. 
\begin{equation*}
\Omega \left( c\right) =\left\{ \dsum_{k}\min \left( 1,\lambda ^{j-k}\right)
c_{k}\right\} _{j}\text{ if }c=\left\{ c_{k}\right\} _{k}.
\end{equation*}%
Then follows readily (see \cite{6-CwPe}, p.35-36, th. 4.8 and th. 4.9) that
we have%
\begin{equation}
\left\Vert a\right\Vert _{A}=\inf_{a}\Psi \left( \left\{ J\left( \lambda
^{k},a\right) \right\} \right) ,  \tag{2}
\end{equation}%
where the inf goes over all representation $a=\left\{ a_{k}\right\} $ of $%
a,\Psi $ being the sequence norm defined as $\Psi =\Phi \circ \Omega $. (All
interpolation spaces are classical J spaces.)$\medskip $

Consider now the special case when $1^{\circ }$ $\overline{A}=l^{\infty
}\left( \overline{w}\right) $ or $2^{\circ }$ $\overline{A}=l^{1}\left( 
\overline{w}\right) $ fixing first attention to the $\lambda $ -adic case,
i.e. $\overline{w}=\left( \left\{ 1\right\} ,\left\{ \lambda ^{-k}\right\}
\right) $ (so that $\overline{A}=\overline{l}_{\lambda }^{\infty }$ or $%
\overline{A}=\overline{l}_{\lambda }^{1}$ respectively in the notation of
sec. 3). $\medskip $

It is easy to see that in both cases we have a $\sigma $-additive Calderon
couple and that in the second case also the strong form of the fundamental
lemma holds. \footnote{%
See notes Appendix: $\left\langle 2\right\rangle .$}.$\medskip $

In case $1^{\circ },$i.e. $\overline{A}=\overline{l}_{\lambda }^{\infty },$
substituting $a$ for the sequence $a^{^{\prime }}=\left\{ K\left( \lambda
^{j},a\right) \right\} $ in $\left( 1\right) ,$we then get since $K\left(
\lambda ^{j},a^{^{\prime }}\right) \approx K\left( \lambda ^{j},a\right) :$%
\begin{equation*}
\left\Vert \left\{ K\left( \lambda ^{j},a\right) \right\} \right\Vert
_{A}\approx \Phi \left( \left\{ K\left( \lambda ^{j},a\right) \right\}
\right) .
\end{equation*}%
so using $\left( 1\right) $ once more we see that%
\begin{equation}
\left\Vert a\right\Vert _{A}\approx \left\Vert \left\{ K\left( \lambda
^{j},a\right) \right\} \right\Vert _{A}.  \tag{3}
\end{equation}%
In the same way in case $2^{\circ },$i.e. $\overline{A}=\overline{l}%
_{\lambda }^{1},$ we obtain using $\left( 2\right) $%
\begin{equation}
\left\Vert a\right\Vert _{A}\approx \inf_{a}\left\Vert \left\{ J\left(
\lambda ^{j},a_{j}\right) \right\} \right\Vert _{A}.  \tag{4}
\end{equation}%
It is easy to see that upon introducing an equivalent norm in $A$ we can
achieve equality in $\left( 3\right) $ or $\left( 4\right) .$ (That much is
of course not needed for the aims we are attempting at.)$\medskip $

We have thus established \underline{1}, formulas $\left( 3\right) $ and $%
\left( 4\right) $ in the special case of $\lambda $-adic sequences.$\medskip 
$

The general case can be handled similarly by observing first that the
K-functional $K\left( t,a\right) $ is up to equivalence determined by its
restriction to the sequence $\left\{ \tau _{n}\right\} .$(Recall that $\tau
_{n}=w_{n}^{0}/w_{n}^{1}.).\medskip $

There remains to establish \underline{1}, formulas $\left( 5\right) $ and $%
\left( 6\right) .$ With no loss of generality we may take $%
w_{n}^{1}=1,w_{n}^{1}=1/\tau _{n}.$ As for \underline{1}, $\left( 5\right) $
we notice as is well-known (see e.g. \cite{30-Pe}), that%
\begin{equation*}
K\left( t,a\right) \approx \sup_{n}\min \left( 1,\frac{t}{\tau _{n}}\right)
\left\vert a_{n}\right\vert \text{ if }a=\left\{ a_{n}\right\} .
\end{equation*}%
We then see that%
\begin{eqnarray*}
\sup_{n}\frac{K\left( \tau _{n},a\right) }{\rho \left( \tau _{n}\right) }
&\approx &\sup_{m}\sup_{n}\frac{\min \left( 1,\frac{\tau _{m}}{\tau _{n}}%
\right) }{\rho \left( \tau _{n}\right) }\left\vert a_{m}\right\vert  \\
&\approx &\sup_{m}\frac{\left\vert a_{m}\right\vert }{\rho \left( \tau
_{m}\right) }=\left\Vert a\right\Vert _{l^{\infty }\left( \rho \left( 
\overline{w}\right) \right) }.
\end{eqnarray*}%
Here we have effectively made use of $\rho \in \mathcal{P}$. This proves 
\underline{1}, $\left( 5\right) .$ \underline{1},$\left( 6\right) $ can be
proven along similar lines. we omit the details (cf. \cite{6-CwPe}).

\section*{Notes}

\underline{Introduction.}

$\left\langle 1\right\rangle :$ After most of this was completed there
appeared the note by Brudnyi and Krugjlak \cite{4-BrKr} listing a number of
even more spectacular results. In particular the "strong form of the
fundamental lemma" (cf. e.g. the appendix of the present paper) seems to be
applicable to all couples with practically no reserve. As a result large
portions of this paper could have been presented in a neater form. See also
the forthcoming paper by Nilsson \cite{22-Nil} where the results of Brudnyi
and Krugjlak are further exploited.

$\left\langle 2\right\rangle .$ The same result appears also in the note by
Brudnyi and Krugjlak \cite{4-BrKr}.$\medskip $

\underline{$\frac{1}{2}$} .

$\left\langle 1\right\rangle :$ If $\overline{X}$ is regular then $\left(
\Sigma \overline{X}\right) ^{^{\prime }}\approx \Delta \overline{X}%
^{^{\prime }}$ where $\overline{X}^{^{\prime }}=\left( X_{0}^{^{\prime
}},X_{1}^{^{\prime }}\right) $ is the dual pair. If $\overline{X}$ is not
regular then $\overline{X}^{^{\prime }}$ is not formally defined but we have
at any rate always canonical maps $\left( \Sigma \overline{X}\right)
^{^{\prime }}\rightarrow X_{j}^{^{\prime }},\left( j=0,1\right) $ which
however need not be injective; $\left\Vert l_{n}\right\Vert
_{X_{j}^{^{\prime }}}$ in $\left( 1\right) $ is thus an abusive notation for
the norm of the of the image of $l_{n}.$

$\left\langle 2\right\rangle .$ By abuse of notation we now write $T$ in
place of $T_{i}=T\mid _{X_{i}}.$We intend to follows this praxis in what
follows.

$\left\langle 3\right\rangle .$ Following \cite{6-CwPe} we say that we have
a representation of $x$ (where $x\in \Sigma _{0}\overline{X})$ if there is
given a sequence $\widehat{x}=\left\{ x_{v}\right\} _{v\in \mathbb{\mathbb{Z}%
}}$ in $\Delta \overline{X}$ such that $x=\dsum x_{v}$ (convergence in $%
\Sigma \overline{X}).$

$\left\langle 4\right\rangle :$ For any sequence $w=\left\{ w_{v}\right\}
_{v\in \mathbb{Z}}$ of positive real numbers and $p\in \left[ 1,\infty %
\right] $ we let $l^{p}\left( w\right) $ be the set of real sequences $%
a=\left\{ a_{k}\right\} _{k\in \mathbb{Z}}$such that $\left( \dsum_{k\in 
\mathbb{\mathbb{Z}}}\left( w_{k}\left\vert a_{k}\right\vert \right)
^{p}\right) ^{1/p}<\infty $, with the usual interpretation if $p=\infty
.\medskip $

\underline{2}.

$\left\langle 1\right\rangle .$ This really gives an indication that one
should have worked with the continuous version, not the discrete one.

$\left\langle 2\right\rangle .$ This notation is inspired by Ovchinnikov 
\cite{25-Ovc}.

$\left\langle 3\right\rangle .$I am grateful to Lars G\aa rding for
directing my attention to Godement's work.$\medskip $

\underline{3}.

$\left\langle 1\right\rangle .$ If $\left( 1\right) $ and $\left( 2\right) $
are not stipulated one usual takes the former condition as a definition; cf. 
\cite{17-KrPeSe}.$\medskip $

\underline{4}.

$\left\langle 1\right\rangle :$ Cwikel \cite{5-Cw} begins also the
investigation of those couples $\overline{X}$ and $\overline{Y}$ for which
the same conclusion holds if we substitute $\dsum_{k\in \mathbb{\mathbb{Z}}%
}c_{k}<\infty $ for $\left( \dsum_{k\in \mathbb{Z}}c_{k}^{r}\right)
^{1/r}<\infty ,r$ a fixed number $>1.$(If we can pass to the limit $r=\infty 
$ we have a Calderon pair.). The infimum of the numbers $r$ which are
permissible here apparently is an important invariant for the couples $%
\overline{X}$ and $\overline{Y}.$

$\left\langle 2\right\rangle .$ Ovchinnikov quotes also the paper by
Lindenstrauss-Pelczynski \cite{18-LiPe} - the main \underline{raison d'etre}
of the latter was indeed to make Grothendieck's work available for a general
audience - but as far as we have been able to see this particular
consequence is not, at least not very explicitly, treated there.

$\left\langle 3\right\rangle .$ The nuclearity. of $R:B\rightarrow
l^{1}\left( w\right) ,$ where $B$ is any Banach space implies indeed that $%
\dsum_{k}\left\Vert ^{t}R\left( e_{k}\right) \right\Vert _{B^{^{\prime
}}}w_{k}<\infty $ where $e_{k}$ is the k'th basis vector of $l^{1}\left(
w\right) .\medskip $

\underline{5}.

$\left\langle 1\right\rangle .$ Cf. \cite{19-LiTz}, prop. 1.c.8 where a
similar device is used; I owe this observation to Per Nilsson.

$\left\langle 2\right\rangle .$We borrow this term from \cite{17-KrPeSe}. $%
\medskip $

\underline{Appendix}.

$\left\langle 1\right\rangle .$From the work of Brudnyi-Krugljak. \cite%
{4-BrKr} referred to already in the Introduction we infer that the
assumption of $\sigma $ - additivity is essentially superfluous.

$\left\langle 2\right\rangle .$That we have a Calderon couple results from 
\cite{27-Pe} (cf. \cite{28-Pe}, \cite{6-CwPe} and \cite{31-SeSe}
respectively). The rest is easy to verify.$\medskip $

\end{document}